\documentclass[journal]{IEEEtran}
\IEEEoverridecommandlockouts
\usepackage[cmex10,intlimits]{amsmath}
\usepackage{epsf}
\usepackage[pdftex]{graphicx}
\usepackage{tikz,pgfplots,tikz-3dplot}
\usetikzlibrary{shapes.geometric}
\usetikzlibrary{shapes.geometric, arrows}
\usetikzlibrary{decorations.shapes}
\usepackage{bigstrut}
\usepackage[noadjust]{cite}
\usepackage{mathtools}
\usepackage{algorithm}
\usepackage{algorithmic}
\usepackage{amssymb}
\usepackage{enumitem}
\usepackage{esvect}

\usepackage{multirow}
\usepackage{multicol}
\usepackage{array}
\newcolumntype{L}[1]{>{\raggedright\let\newline\\\arraybackslash\hspace{0pt}}m{#1}}
\newcolumntype{C}[1]{>{\centering\let\newline\\\arraybackslash\hspace{0pt}}m{#1}}
\newcolumntype{R}[1]{>{\raggedleft\let\newline\\\arraybackslash\hspace{0pt}}m{#1}}

\usepackage{amsfonts}
\usepackage{epsfig}
\usepackage{textcomp}
\usepackage{mleftright}

\newcommand{\gv}[1]{\ensuremath{\mbox{\boldmath$ #1 $}}}

\newcommand{\grad}[1]{\gv{\nabla} #1}
\newcommand{\abs}[1]{\left| #1 \right|}

\allowdisplaybreaks

\def\BibTeX{{\rm B\kern-.05em{\sc i\kern-.025em b}\kern-.08em
    T\kern-.1667em\lower.7ex\hbox{E}\kern-.125emX}}

\usepackage[subrefformat=parens,labelformat=parens,caption=false]{subfig}
\begin{document}

\title{An Efficient Integral Equation Method for Full-wave Analysis of Inhomogeneous Electromagnetic Surfaces with Connected Conductors}

\author{Reza Gholami,~\IEEEmembership{Member,~IEEE},~Parinaz Naseri,~\IEEEmembership{Student Member,~IEEE},~Piero Triverio,~\IEEEmembership{Senior Member,~IEEE}, and~Sean V. Hum,~\IEEEmembership{Senior Member,~IEEE}

\thanks{Manuscript received Apr~23, 2021;}
\thanks{The work was supported by a Strategic Partnerships Grant from Natural Sciences and Engineering Research Council (NSERC) and CMC Microsystems.}
\thanks{The authors are with the Department of Electrical and Computer Engineering, University of Toronto, Toronto, Canada (e-mail: r.gholami@utoronto.ca, parinaz.naseri@utoronto.ca, piero.triverio@utoronto.ca, and sean.hum@utoronto.ca).}
\thanks{Color versions of one or more of the figures in this paper are available online at http://ieeexplore.ieee.org.}
\thanks{Digital Object Identifier}
}

\markboth{IEEE Transactions on Antennas and Propagation, Vol. XX, No. XX, April 2021}{Gholami, Triverio, Hum, A New Domain Decomposition Technique for Full-wave Analysis of Inhomogeneous Electromagnetic Surfaces with Connected Regions}
\maketitle

\begin{abstract}
In this paper, a generalized macromodeling approach is presented to simulate complex electromagnetic (EM) surfaces consisting of unit cells with connected conductors. Macromodels of each unit cell are produced by applying the equivalence principle on fictitious surfaces encapsulating them.
Unit cells often consist of multiple dielectric layers and conductor traces, featuring multiscale structures. Challenges arise when a current-carrying conductor trace traverses the fictitious surface.
Hence, a new method based on half Rao-Wilton-Glisson basis functions is proposed to accurately ensure the continuity of the surface currents and avoid singularities at the intersections. 
The accuracy of the proposed approach is validated by comparing the results with commercial solvers for different EM surfaces.

\end{abstract}

\IEEEpeerreviewmaketitle

\section{Introduction}
Full-wave analysis of complex electromagnetic (EM) surfaces, such as frequency selective surfaces, reflectarrays, transmitarrays, and metasurfaces 
has been critically important to the development of wireless communication applications. These surfaces are usually electrically large and are often composed arrangements of finely-textured conductors within thin, multilayer dielectrics
with subwavelength dimensions. Therefore, their analysis is a time-consuming and resource-demanding problem that becomes the bottleneck of their design and optimization. 

Integral equation techniques~\cite{Chew_IE_Book, PMCHWT_Muller_Ergul} can be a powerful tool to accurately simulate such problems. The method of moments (MoM)~\cite{Harrington-MoM, Gibson-MoM} discretization of these large-scale problems results in a dense matrix equation which is computationally very challenging to solve. Employing sophisticated fast iterative~\cite{MLFMM-Chew, Fast Algorithms book, MLFMM-Ergul-Gurel, Bleszynski-AIM, Yilmaz-VIE-AIM} and direct solution~\cite{Hmatrix book_Hackbusch, Shaeffer-algorithm, Boag, ACA_Zhao} algorithms  effectively reduces the amounts of memory and CPU time required while producing accurate results. While iterative fast methods such as the multilevel fast multipole method (MLFMM)~\cite{MLFMM-Ergul-Gurel} and the adaptive integral method (AIM)~\cite{Bleszynski-AIM} are popular for such problems, they can suffer from slow convergence. This is due to the fact that the multiscale features of the EM structures deteriorate the matrix conditioning. While adding appropriate preconditioners can circumvent this problem, it often reduces the efficiency of the EM solver due to the expensive inversion operation. Alternatively, fast direct methods such as the hierarchical ($\mathcal H$)-matrix~\cite{Hmatrix book_Hackbusch} are largely insensitive to the matrix conditioning due to their non-iterative nature. However, their performance declines rapidly for large-scale problems (usually objects with electrical size larger than 10$\lambda$) due to the growth in the rank of the admissible sub-matrices~\cite{Reza-Thesis}.       
  
The equivalence principle algorithm (EPA)~\cite{EPA-Chew1, EPA-Chew2, EPA-Lancellotti, EPA-Taskinen} is another IE technique that can be very efficient for full-wave analysis of large-scale and highly complex EM surfaces. In the EPA, each unit cell of the EM surface is encapsulated by a simple fictitious surface that generates a subdomain of the overall surface. Using Love's equivalence principle and the surface integral equation (SIE), one can define an operator that encodes the EM behavior of the objects inside the unit cell, relating incident and scattered field on the fictitious surface. 
The coupling between the subdomains is computed using a translation operator. This transfers the unknowns on the unit cell, which may have high mesh density, to the unknowns on the fictitious surface that has far fewer. Therefore, the final MoM matrix in the EPA approach has better convergence and can be solved with fewer computational resources.

Another approach for EM surfaces based on the equivalence principle is the macromodeling technique proposed in~\cite{Utkarsh_macromodel_2018, Utkarsh_macromodel_2020}. In this approach, a complex unit cell is also enclosed by a simple fictitious surface. Initial formulations of the method were only applicable to scatterers comprising conductors~\cite{Utkarsh_macromodel_2018}, while later development allowed them to 
handle composite dielectric-conductor scatterers~\cite{Utkarsh_macromodel_2020}. A macromodel operator is generated 
to capture the scattering properties of the unit cell by equivalent electric and magnetic currents introduced on the fictitious surface. The interelement coupling is computed using the electric field integral equation (EFIE) and the magnetic field integral equation (MFIE). 


The main challenge in the EPA approach is when a conductor trace traverses the fictitious surfaces. This is the case for many practical EM surfaces with tailored scattering properties, e.g. EM surfaces realized from screens (complementary surfaces), connected arrays, etc.~\cite{EMS_Connected1, EMS_Connected2, EMS_Connected3, EMS_Connected4, EMS_Connected5, EMS_Connected6, EMS_Connected7, EMS_Connected8, EMS_Connected9}. Breaking such conducting scatterers into two pieces causes a current discontinuity which produces singularity in the fields and consequently produces inaccurate solutions~\cite{EPA-Chew2, EPA-SRM}. The tap basis method was implemented in the EPA to mitigate the singularity of the current between scatterers connected by conductor traces~\cite{EPA-Chew1, EPA-Chew2, TAP1}. While the tap basis method is a sufficient approach to model the continuity of current in EPA, it deteriorates the conditioning of the EPA equations~\cite{EPA-SRM}.


In this work, we generalize the macromodeling approach~\cite{Utkarsh_macromodel_2020} previously developed for the case of EM surfaces with unit cells composed of non-connected conductor scatterers to the case of EM surfaces where the unit cells are connected with conductors. To model the continuity of the current and avoid the singularity in the macromodeling approach, a new technique based on half Rao-Wilton-Glisson (RWG) basis functions is proposed. The new method offers a simple way to handle conductive traces traversing the macromodel boundaries, while preserving the conditioning of the original formulation. In order to enable the analysis of large problems, a previously developed acceleration technique based on the fast Fourier transform (FFT) is extended to the case of surfaces with connected conductors. The proposed technique is validated through numerical analysis of multiple scattering problems and comparison with commercial EM solvers. It is shown that the proposed macromodeling approach combined with FFT-based acceleration can efficiently reduce the time and resources needed to analyze large EM surfaces for various applications.

The rest of paper is organized as follows. In Section~\ref{Overview of Macromodeling Approach}, we review the previously developed macromodeling approach and discuss its limitation. In Section~\ref{proposed method}, we describe the proposed technique. Numerical results and discussion in Section~\ref{Numerical Results} demonstrate the validity and efficiency of the proposed method. Finally, a summary and conclusion are included
in Section~\ref{Conclusion}.

\section{Overview of Macromodeling Approach}
\label{Overview of Macromodeling Approach}
In this section, we review the different steps in the macromodeling technique developed previously~\cite{Utkarsh_macromodel_2020} for the analysis of EM surfaces. 

\subsection{Surface Integral Equations}
\label{Surface Integral Equation}
\begin{figure}[t]
\hfill \null
\\
\null \hfill
\subfloat[ \label{fig:sample_cell_equivalent_1_oldpaper} Original problem] {\resizebox{0.38\columnwidth} {!} {\input{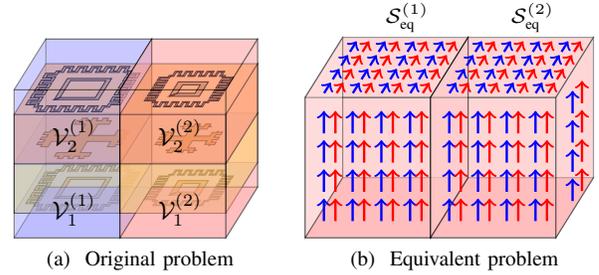}}}
\hfill 
\subfloat[\label{fig:sample_cell_equivalent_2_oldpaper} Equivalent problem]{\resizebox{0.45\columnwidth} {!}  {\begin{tikzpicture}[xscale= 0.25, yscale =0.4]
\pgfmathsetmacro{\h}{2*0.7}

\pgfmathsetmacro{\UnitcellWidth}{3.75}
\pgfmathsetmacro{\UnitcellLength}{3.75}
\pgfmathsetmacro{\UnitcellHeight}{3.75}
\pgfmathsetmacro{\numarrow}{3}
\pgfmathsetmacro{\numarrowa}{4}


\begin{scope}[shift = {($(0,-1.5*\h)$)}, every node/.append style={xslant=1},xslant=1, xscale= 1, yscale = 0.3]
\coordinate (p1) at (-3.0,-3.0);
\coordinate (p2) at (3.0,-3.0);
\coordinate (p3) at (3.0, 3.0);
\coordinate (p4) at (-3.0,3.0);
\draw[fill=red!40!white,opacity=0.4] (p1)--(p2)--(p3)--(p4)--(p1);
\end{scope}

\draw[fill=red!40!white, opacity = 0.4] (p1) --++ (0,3*\h) --++ (6.0,0) -- (p2) -- (p1);

\begin{scope}[shift = {($(0,1.5*\h)$)}, every node/.append style={xslant=1},xslant=1, xscale= 1, yscale = 0.3]
\coordinate (p1) at (-3.0,-3.0);
\coordinate (p2) at (3.0,-3.0);
\coordinate (p3) at (3.0, 3.0);
\coordinate (p4) at (-3.0,3.0);
\draw[fill=red!40!white,opacity=0.4] (p1)--(p2)--(p3)--(p4)--(p1);

\end{scope}

\draw[fill=red!40!white, opacity = 0.4] (p2) -- (p3) -- ($(p3) + (0,-3*\h)$) -- ($(p2) + (0,-3*\h)$) -- (p2);

\draw[fill=red!10!white, opacity = 0.4] (p1) -- (p4) -- ($(p4) + (0,-3*\h)$) -- ($(p1) + (0,-3*\h)$) -- (p1);

\begin{scope}[shift={(-1.5,-1.0)}]
\foreach \x in {0,...,\numarrow}
{
\foreach \y in {0,...,\numarrow}
{
	\coordinate (t1) at (-\UnitcellWidth*0.4 + \x*1.1*\UnitcellWidth/\numarrow, -\UnitcellLength*0.4 + \y*0.7*\UnitcellLength/\numarrow);
	\coordinate (t2) at (-\UnitcellWidth*0.4 + \x*1.1*\UnitcellWidth/\numarrow, -\UnitcellLength*0.4 + \y*0.7*\UnitcellLength/\numarrow + 0.18*\UnitcellLength);
\draw [blue, ->, thick] (t1) -- (t2);
}
}
\begin{scope}[shift={(0.55,0.0)}]
\foreach \x in {0,...,\numarrow}
{
\foreach \y in {0,...,\numarrow}
{
	\coordinate (t1) at (-\UnitcellWidth*0.4 + \x*1.1*\UnitcellWidth/\numarrow, -\UnitcellLength*0.4 + \y*0.7*\UnitcellLength/\numarrow);
	\coordinate (t2) at (-\UnitcellWidth*0.4 + \x*1.1*\UnitcellWidth/\numarrow, -\UnitcellLength*0.4 + \y*0.7*\UnitcellLength/\numarrow + 0.18*\UnitcellLength);
\draw [red, ->, thick] (t1) -- (t2);
}
}
\end{scope}

\end{scope}

\begin{scope}[shift = {($(-1.3,0.8*\h)$)}, every node/.append style={xslant=1},xslant=1, xscale= 1, yscale = 0.3]
\begin{scope}[shift={(-0.4,2.2)}]
\foreach \x in {0,...,\numarrow}
{
\foreach \y in {0,...,\numarrow}
{
	\coordinate (t1) at (-\UnitcellWidth*0.4 + \x*1.1*\UnitcellWidth/\numarrow, -\UnitcellLength*0.4 + \y*1.1*\UnitcellLength/\numarrow);
	\coordinate (t2) at (-\UnitcellWidth*0.4 + \x*1.1*\UnitcellWidth/\numarrow, -\UnitcellLength*0.4 + \y*1.1*\UnitcellLength/\numarrow + 0.28*\UnitcellLength);
\draw [blue, ->, thick] (t1) -- (t2);
}
}
\end{scope}

\begin{scope}[shift={(0.2,2.2)}]
\foreach \x in {0,...,\numarrow}
{
\foreach \y in {0,...,\numarrow}
{
	\coordinate (t1) at (-\UnitcellWidth*0.4 + \x*1.1*\UnitcellWidth/\numarrow, -\UnitcellLength*0.4 + \y*1.1*\UnitcellLength/\numarrow);
	\coordinate (t2) at (-\UnitcellWidth*0.4 + \x*1.1*\UnitcellWidth/\numarrow, -\UnitcellLength*0.4 + \y*1.1*\UnitcellLength/\numarrow + 0.28*\UnitcellLength);
\draw [red, ->, thick] (t1) -- (t2);
}
}
\end{scope}

\end{scope}

\begin{scope}[shift={(6.0,0.0)}]
\begin{scope}[shift = {($(0,-1.5*\h)$)}, every node/.append style={xslant=1},xslant=1, xscale= 1, yscale = 0.3]
\coordinate (p1) at (-3.0,-3.0);
\coordinate (p2) at (3.0,-3.0);
\coordinate (p3) at (3.0, 3.0);
\coordinate (p4) at (-3.0,3.0);
\draw[fill=red!40!white,opacity=0.4] (p1)--(p2)--(p3)--(p4)--(p1);
\end{scope}

\draw[fill=red!40!white, opacity = 0.4] (p1) --++ (0,3*\h) --++ (6.0,0) -- (p2) -- (p1);

\begin{scope}[shift = {($(0,1.5*\h)$)}, every node/.append style={xslant=1},xslant=1, xscale= 1, yscale = 0.3]
\coordinate (p1) at (-3.0,-3.0);
\coordinate (p2) at (3.0,-3.0);
\coordinate (p3) at (3.0, 3.0);
\coordinate (p4) at (-3.0,3.0);
\draw[fill=red!40!white,opacity=0.4] (p1)--(p2)--(p3)--(p4)--(p1);
\end{scope}

\draw[fill=red!40!white, opacity = 0.4] (p2) -- (p3) -- ($(p3) + (0,-3*\h)$) -- ($(p2) + (0,-3*\h)$) -- (p2);

\draw[fill=red!10!white, opacity = 0.4] (p1) -- (p4) -- ($(p4) + (0,-3*\h)$) -- ($(p1) + (0,-3*\h)$) -- (p1);

\begin{scope}[shift={(-1.5,-1.0)}]
\foreach \x in {0,...,\numarrow}
{
\foreach \y in {0,...,\numarrow}
{
	\coordinate (t1) at (-\UnitcellWidth*0.4 + \x*1.1*\UnitcellWidth/\numarrow, -\UnitcellLength*0.4 + \y*0.7*\UnitcellLength/\numarrow);
	\coordinate (t2) at (-\UnitcellWidth*0.4 + \x*1.1*\UnitcellWidth/\numarrow, -\UnitcellLength*0.4 + \y*0.7*\UnitcellLength/\numarrow + 0.18*\UnitcellLength);
\draw [blue, ->, thick] (t1) -- (t2);
}
}

\begin{scope}[shift={(0.55,0.0)}]
\foreach \x in {0,...,\numarrow}
{
\foreach \y in {0,...,\numarrow}
{
	\coordinate (t1) at (-\UnitcellWidth*0.4 + \x*1.1*\UnitcellWidth/\numarrow, -\UnitcellLength*0.4 + \y*0.7*\UnitcellLength/\numarrow);
	\coordinate (t2) at (-\UnitcellWidth*0.4 + \x*1.1*\UnitcellWidth/\numarrow, -\UnitcellLength*0.4 + \y*0.7*\UnitcellLength/\numarrow + 0.18*\UnitcellLength);
\draw [red, ->, thick] (t1) -- (t2);
}
}
\end{scope}

\end{scope}

\begin{scope}[shift = {($(-1.3,0.8*\h)$)}, every node/.append style={xslant=1},xslant=1, xscale= 1, yscale = 0.3]
\begin{scope}[shift={(-0.4,2.2)}]
\foreach \x in {0,...,\numarrow}
{
\foreach \y in {0,...,\numarrow}
{
	\coordinate (t1) at (-\UnitcellWidth*0.4 + \x*1.1*\UnitcellWidth/\numarrow, -\UnitcellLength*0.4 + \y*1.1*\UnitcellLength/\numarrow);
	\coordinate (t2) at (-\UnitcellWidth*0.4 + \x*1.1*\UnitcellWidth/\numarrow, -\UnitcellLength*0.4 + \y*1.1*\UnitcellLength/\numarrow + 0.28*\UnitcellLength);
\draw [blue, ->, thick] (t1) -- (t2);
}
}
\end{scope}

\begin{scope}[shift={(0.2,2.2)}]
\foreach \x in {0,...,\numarrow}
{
\foreach \y in {0,...,\numarrow}
{
	\coordinate (t1) at (-\UnitcellWidth*0.4 + \x*1.1*\UnitcellWidth/\numarrow, -\UnitcellLength*0.4 + \y*1.1*\UnitcellLength/\numarrow);
	\coordinate (t2) at (-\UnitcellWidth*0.4 + \x*1.1*\UnitcellWidth/\numarrow, -\UnitcellLength*0.4 + \y*1.1*\UnitcellLength/\numarrow + 0.28*\UnitcellLength);
\draw [red, ->, thick] (t1) -- (t2);
}
}
\end{scope}

\end{scope}

\begin{scope}[shift={(4.4,-0.4)}]
\foreach \x in {0,...,0}
{
\foreach \y in {0,...,\numarrow}
{
	\coordinate (t1) at (-\UnitcellWidth*0.4 + \x*0.8*\UnitcellWidth/\numarrow, -\UnitcellLength*0.4 + \y*0.7*\UnitcellLength/\numarrow);
	\coordinate (t2) at (-\UnitcellWidth*0.4 + \x*0.8*\UnitcellWidth/\numarrow, -\UnitcellLength*0.4 + \y*0.7*\UnitcellLength/\numarrow + 0.18*\UnitcellLength);
\draw [blue, ->, thick] (t1) -- (t2);
}
}
\end{scope}

\begin{scope}[shift={(4.85,-0.1)}]
\foreach \x in {0,...,0}
{
\foreach \y in {0,...,\numarrow}
{
	\coordinate (t1) at (-\UnitcellWidth*0.4 + \x*0.8*\UnitcellWidth/\numarrow, -\UnitcellLength*0.4 + \y*0.7*\UnitcellLength/\numarrow);
	\coordinate (t2) at (-\UnitcellWidth*0.4 + \x*0.8*\UnitcellWidth/\numarrow, -\UnitcellLength*0.4 + \y*0.7*\UnitcellLength/\numarrow + 0.18*\UnitcellLength);
\draw [red, ->, thick] (t1) -- (t2);
}
}
\end{scope}

\end{scope}

\draw node at (1.0,3.6) {\scriptsize $\mathcal S_{\text{eq}}^{(1)}$};
\draw node at (7.0,3.6) {\scriptsize $\mathcal S_{\text{eq}}^{(2)}$};

\end{tikzpicture}
}}
\hfill \null
\caption{(a) An example of an array composed of two unit cells simulated by the previous macromodeling approach~\cite{Utkarsh_macromodel_2020}. Each unit cell has two dielectric regions $\mathcal V_1^{(m)}$ and $\mathcal V_2^{(m)}$, $m=1,2$, (shown with yellow and red colors) and is enclosed by a fictitious surface. The conductor traces are not allowed to traverse the fictitious surfaces. (b) Unit cells are modeled by equivalent surface electric and magnetic current densities (shown with blue and red arrows, respectively) introduced on $\mathcal S_{\text{eq}}^{(1)}$ and $\mathcal S_{\text{eq}}^{(2)}$.}
\label{fig:equivalent_cell_oldpaper}
\end{figure}

The macromodeling approach for the analysis of an EM surface starts with generating a macromodel describing the EM response of each unit cell. Let us consider a surface consisting of $M$ unit cells, such as the one shown in Fig.~\ref{fig:equivalent_cell_oldpaper}\subref{fig:sample_cell_equivalent_1_oldpaper}, each with $V$ dielectric layers and an arbitrary number of conductive scatterers. In Fig.~\ref{fig:equivalent_cell_oldpaper}\subref{fig:sample_cell_equivalent_1_oldpaper}, $M=2$ and $V=2$. We first enclose the $m$-th unit cell with a fictitious surface $\mathcal S_\text{eq}^{(m)}$ such that all conductive scatterers and dielectric layers are inside the fictitious surface. Using the equivalence principle~\cite{Chew_IE_Book}, the field inside or outside a homogeneous region within the unit cell, such as $\mathcal V_1^{(m)}$ or $\mathcal V_2^{(m)}$, can be computed through the equivalent surface electric current density $\vec{J}_{v}$ and magnetic current density $\vec{M}_{v}$ defined on the boundary of the $v$-th region. Hence, we discretize all boundaries with triangular elements and expand the surface current densities with RWG basis functions. By applying the EFIE and the MFIE~\cite{Gibson-MoM} for each region in a given unit cell and testing the result with RWG functions, we obtain the system of equations

\begin{equation}
\label{SIE Matrix equation-1}
\begin{aligned}
\underbrace{ \begin{bmatrix} 
\begin{pmatrix} \textbf L_1^E  & \mspace{-8mu}\textbf K_1^E \\ \textbf K_1^H  & \mspace{-8mu}\textbf L_1^H \end{pmatrix} & \mspace{-6mu}0 & \mspace{-6mu}0 \\ 
 0\mspace{-6mu} &  \ddots  & \mspace{-6mu}0 \\ 
0 & \mspace{-6mu}0 & \mspace{-6mu}\begin{pmatrix} \textbf L_V^E  & \mspace{-8mu}\textbf K_V^E \\ \textbf K_V^H  & \mspace{-8mu}\textbf L_V^H \end{pmatrix}
\end{bmatrix}\mspace{-3mu}}_{{\textbf{Z}^{(m)}}}
\underbrace{ \begin{bmatrix}
\begin{pmatrix} \textbf J_1  \\ \textbf M_1  \end{pmatrix} \\
 \vdots  \\
\begin{pmatrix} \textbf J_V  \\ \textbf M_V  \end{pmatrix}
\end{bmatrix}}_{\textbf{X}^{(m)}} \mspace{-3mu}=\mspace{-3mu} 
\begin{bmatrix}
\begin{pmatrix} 0  \\ 0  \end{pmatrix} \\
 \vdots  \\
\begin{pmatrix} 0  \\ 0  \end{pmatrix}
\end{bmatrix}
\end{aligned}
\end{equation}
where $\textbf J_v$ and $\textbf M_v$ are vectors of unknowns collecting the electric and magnetic current density coefficients on the boundary of the ${v}$-th homogeneous region, respectively. In~\eqref{SIE Matrix equation-1}, $\textbf L^{\alpha}_v$ and $\textbf K^{\alpha}_v$ are matrices obtained from the surface integral operators (${\alpha} = E~\text{or}~H$ represents the EFIE or the MFIE, respectively) 
\begin{equation}
\begin{aligned}
\hat{n} \times \mspace{-3mu} \left[ {\mathcal{\vec L}}_v \vec{X}(\vec{r}\,') \right]\mspace{-3mu}(\vec{r}) \mspace{-3mu}=\mspace{-3mu} \hat{n} \times \left[1+\frac{\grad\grad \cdot}{k^2_v} \right]  \int\limits_{\mathcal S_{v}}{\mspace{-6mu}{G}_{v}}(\vec{r},\vec{r}\,') \vec{X}(\vec{r}\,') d\mathcal V^\prime 
\label{L Operator-1}
\end{aligned}
\end{equation} 
  
\begin{equation}
\begin{aligned}
\hat{n} \times \mspace{-3mu} \left[ {\mathcal{\vec K}}_v \vec{X}(\vec{r}\,') \right]\mspace{-3mu}(\vec{r}) \mspace{-3mu}= \mspace{2mu}& \hat{n} \times \text{p.v.} \left[\grad \mspace{-3mu} \times \mspace{-6mu} \int\limits_{\mathcal S_{v}}{{G}_{v}}(\vec{r},\vec{r}\,') \vec{X}(\vec{r}\,') d\mathcal V^\prime \right] \\
&+\frac{\vec{X}(\vec{r}\,')}{2} 
\label{K Operator-1}
\end{aligned}
\end{equation} 
with $k_v=\omega \sqrt{\varepsilon_v \mu_0}$ being the wavenumber of the $v$-th dielectric region, and $\varepsilon_v$ being its electrical permittivity. In \eqref{L Operator-1} and \eqref{K Operator-1}, ${{G}}_{v}(\vec{r},\vec{r}\,')=\text{exp}(-j k_v \abs{\vec{r}-\vec{r}\,'})/(4\pi\abs{\vec{r}-\vec{r}\,'})$ is the homogeneous Green's function of $v$-th region with $\mathcal S_{v}$ being its boundary. The operator p.v. indicates \textit{principal value}.    

\subsection{Enforcement of Boundary Conditions Inside Each Macromodel} 
The electromagnetic boundary conditions must be imposed on the interface of two regions of a unit cell and at the junctions between the unit cell and its fictitious surface $\mathcal S_{\text{eq}}^{(m)}$. 
The boundary conditions can be applied through~\cite{Utkarsh_macromodel_2020}

\begin{equation}
\begin{aligned}
\textbf  X^{(m)}=\textbf U^{(m)} \widetilde{\textbf{X}}^{(m)},
\label{UX}
\end{aligned}
\end{equation} 
where $\widetilde{\textbf X}^{(m)}=\left[ \widetilde{\textbf X}_{\text{eq}}^{(m)} ~~ \widetilde{\textbf X}_{\text{int}}^{(m)} \right]^T$ is a set of unique unknowns for the $m$-th unit cell that includes both $\widetilde{\textbf X}_{\text{eq}}^{(m)}$ as the unknown current coefficients on $\mathcal S_{\text{eq}}^{(m)}$ and $\widetilde{\textbf X}_{\text{int}}^{(m)}$ as the unknown current coefficients on the surfaces inside $\mathcal S_{\text{eq}}^{(m)}$. Also, $\textbf U^{(m)}$ is a sparse matrix with a few entries per row with values $+1$ or $-1$~\cite{Utkarsh_macromodel_2020}. Enforcement of the boundary conditions will remove the extra unknowns on the interfaces between two or more regions and between the unit cell and the fictitious surface arranged inside the vector of original unknowns $\textbf X^{(m)}$ in~\eqref{SIE Matrix equation-1}. As a result, upon the substitution of~\eqref{UX} into~\eqref{SIE Matrix equation-1} following the left-multiplication by ${\textbf U^{(m)}}^{T}$, we obtain

\begin{equation}
\label{UtZUX_matrix-1}
\begin{aligned}
{\textbf U^{(m)}}^T \textbf Z^{(m)} \textbf U^{(m)} \widetilde{\textbf X}^{(m)} = \begin{bmatrix} \textbf Z_{\text{eq},\text{eq}}^{(m)} & \textbf Z_{\text{eq},\text{int}}^{(m)}  \\ \textbf Z_{\text{int},\text{eq}}^{(m)} & \textbf Z_{\text{int},\text{int}}^{(m)}  \end{bmatrix} \begin{bmatrix} \widetilde{\textbf X}_\text{eq}^{(m)} \\ \widetilde{\textbf X}_\text{int}^{(m)} \end{bmatrix} = \begin{bmatrix} 0 \\ 0 \end{bmatrix}
\end{aligned}.
\end{equation} 

\subsection{Schur Complement}
The next step is to eliminate the interior unknowns $\widetilde{\textbf X}_\text{int}^{(m)}$ in~\eqref{UtZUX_matrix-1} by applying the Schur complement so that  
\begin{equation}
\label{Shur_complement}
\begin{aligned}
\left[ \textbf Z_{\text{eq},\text{eq}}^{(m)} - \textbf Z_{\text{eq},\text{int}}^{(m)} {\textbf Z^{(m)}_{\text{int},\text{int}}}^{-1} {\textbf Z}_{\text{int},\text{eq}}^{(m)} \right] \widetilde{\textbf X}_\text{eq}^{(m)} = \widetilde{\textbf Z}_{\text{eq},\text{eq}}^{(m)}~\widetilde{\textbf X}_\text{eq}^{(m)} = 0.
\end{aligned}
\end{equation} 
The key advantage offered by \eqref{Shur_complement} is the fewer unknowns that are associated with the equivalent surface $S_{\text{eq}}^{(m)}$. This relatively small  system can be computed directly in order to obtain matrix $\widetilde{\textbf Z}_{\text{eq},\text{eq}}^{(m)}$ as the macromodel matrix of the unit cell. The macromodel matrix describes the behavior of the unit cell using only equivalent electric and magnetic surface currents defined on $S_{\text{eq}}^{(m)}$. The above step to generate macromodels only needs to be performed for unique unit cells, which is a significant advantage in analyzing EM surfaces that are composed of a set of unique unit cells~\cite{Utkarsh_macromodel_2020}.


\subsection{Interelement Coupling} 
In order to simulate an EM surface consisting of $M$ unit cells we first need to create the overall macromodel matrix as
\begin{equation}
\label{Overall Macromodel Matrix-1}
\begin{aligned}
\underbrace{ \begin{bmatrix} 
\widetilde{\textbf Z}^{(1)}_{\text{eq},\text{eq}} &  &  \\ 
  &  \ddots  &  \\ 
 &  & \widetilde{\textbf Z}^{(M)}_{\text{eq},\text{eq}}
\end{bmatrix}\mspace{-3mu}}_{{\textbf Z_{\text{eq}}}}
\underbrace{ \begin{bmatrix}
\widetilde{\textbf X}^{(1)}_\text{eq} \\
 \vdots  \\
\widetilde{\textbf X}^{(M)}_\text{eq}
\end{bmatrix}}_{\textbf Y} \mspace{-3mu}=\mspace{-3mu} 
\begin{bmatrix}
0 \\
 \vdots  \\
0
\end{bmatrix}
\end{aligned}
\end{equation} 
where $\textbf Z_{\text{eq}}$ is a block-diagonal matrix known as the overall macromodel matrix and $\textbf Y$ is the overall vector of unknowns. 

So far the unit cells are modeled by the equivalent current densities $\vec J_{\text{eq}}^{~(m)}$ and $\vec M_{\text{eq}}^{(m)}$ defined on $S_{\text{eq}}^{(m)}$. The effect of the fields created from one equivalent surface on the other equivalent surfaces is computed by the EFIE and MFIE for the exterior domain as~\cite{Utkarsh_macromodel_2020}
\begin{equation}
\label{SLAE_exterior EFIE-MFIE-1}
\begin{aligned}
\underbrace{ \begin{bmatrix} 
{\textbf Z}_0^{(1,1)} & {\textbf Z}_0^{(1,2)} & \dots & {\textbf Z}_0^{(1,M)} \\ 
{\textbf Z}_0^{(2,1)} & {\textbf Z}_0^{(2,2)} & \dots & {\textbf Z}_0^{(2,M)} \\ 
\vdots & \vdots & \dots & \vdots \\ 
{\textbf Z}_0^{(M,1)} & {\textbf Z}_0^{(M,2)} & \dots & {\textbf Z}_0^{(M,M)} 
\end{bmatrix}\mspace{-3mu}}_{{\textbf Z_{0}}}
\underbrace{ \begin{bmatrix}
\widetilde{\textbf X}^{(1)}_\text{eq} \\
\widetilde{\textbf X}^{(2)}_\text{eq} \\
 \vdots  \\
\widetilde{\textbf X}^{(M)}_\text{eq}
\end{bmatrix}}_{\textbf Y} \mspace{-3mu}=\mspace{-3mu} 
\underbrace{\begin{bmatrix}
\textbf V_{\phantom{eq}}^{(1)} \\
\textbf V_{\phantom{eq}}^{(2)} \\
 \vdots  \\
\textbf V_{\phantom{eq}}^{(M)} \\
\end{bmatrix}}_{\textbf V} 
\end{aligned},
\end{equation} 
where matrix ${\textbf Z}^{(m,m')}_{0}$ computes the full-wave interaction between the $m$-th and $m'$-th unit cells. Also, $\textbf V^{(m)}$ is the excitation vector which is the projection of the incident electric and magnetic fields on $S_{\text{eq}}^{(m)}$.

\subsection{Enforcement of Boundary Conditions on Equivalent Surfaces}
To produce a well-conditioned matrix equation, we need to couple the overall macromodel matrix equation in \eqref{Overall Macromodel Matrix-1} and interelement coupling matrix equation in \eqref{SLAE_exterior EFIE-MFIE-1} as~\cite{Utkarsh_macromodel_2020}

\begin{equation}
\label{Well-condition formulation-1}
\begin{aligned}
\left(\textbf Z_{\text{eq}} + \textbf Z_{0}\right) \textbf Y = \textbf V.
\end{aligned}
\end{equation} 
Boundary conditions must be imposed at the interface between two or more adjacent equivalent surfaces. Hence, a new sparse matrix $\textbf U_0$ is introduced to remove the redundant unknowns by explicitly enforcing the continuity of tangential electric and magnetic fields on the interfaces between the adjacent equivalent surfaces as 
\begin{equation}
\label{U0Y}
\begin{aligned}
\textbf Y =\textbf  U_{0} \widetilde{\textbf Y},
\end{aligned}
\end{equation}
where $\widetilde{\textbf Y}$ is a vector of unique unknowns. Upon the substitution of~\eqref{U0Y} into~\eqref{Well-condition formulation-1} followed by the left-multiplication of the resulting equations by $\textbf  U_{0}^T$, we obtain~\cite{Utkarsh_macromodel_2020}

\begin{equation}
\label{U0tZU0X_matrix-1}
\begin{aligned}
\textbf U_0^T \left[ \textbf Z_{\text{eq}}+\textbf Z_0 \right] \textbf U_0 \widetilde{\textbf Y} = \underbrace{\textbf U_0^T \textbf V}_{\widetilde{\textbf V}},
\end{aligned}
\end{equation} 
where $\widetilde{\textbf V}$ is the excitation vector obtained after enforcing boundary conditions. The final system of equations in~\eqref{U0tZU0X_matrix-1} can be solved directly or iteratively to obtain current densities on the equivalent surfaces. In the macromodeling approach~\cite{Utkarsh_macromodel_2020}, equation~\eqref{U0tZU0X_matrix-1} is solved iteratively with the generalized minimal residual method (GMRES)~\cite{GMRES1, GMRES2} and accelerated by applying the FFT~\cite{FFT1, FFT2}.


\begin{figure}[t]
\hfill \null
\\
\null \hfill

\subfloat[ \label{fig:sample_cell_equivalent_1} Original problem] {\resizebox{0.45\columnwidth} {!} {\begin{tikzpicture}[xscale= 0.25, yscale =0.4]
\pgfmathsetmacro{\h}{0.7}


\begin{scope}[shift = {($(0,-1.5*\h)$)}, every node/.append style={xslant=1},xslant=1, xscale= 1, yscale = 0.3]
\coordinate (p1) at (-3.0,-3.0);
\coordinate (p2) at (3.0,-3.0);
\coordinate (p3) at (3.0, 3.0);
\coordinate (p4) at (-3.0,3.0);
\draw[fill=red!40!white,opacity=0.4] (p1)--(p2)--(p3)--(p4)--(p1);
\end{scope}

\draw[fill=red!40!white, opacity = 0.4] (p1) --++ (0,3*\h) --++ (6.0,0) -- (p2) -- (p1);

\begin{scope}[shift = {($(0,1.5*\h)$)}, every node/.append style={xslant=1},xslant=1, xscale= 1, yscale = 0.3]
\coordinate (p1) at (-3.0,-3.0);
\coordinate (p2) at (3.0,-3.0);
\coordinate (p3) at (3.0, 3.0);
\coordinate (p4) at (-3.0,3.0);
\draw[fill=red!40!white,opacity=0.4] (p1)--(p2)--(p3)--(p4)--(p1);
\end{scope}

\draw[fill=red!40!white, opacity = 0.4] (p2) -- (p3) -- ($(p3) + (0,-3*\h)$) -- ($(p2) + (0,-3*\h)$) -- (p2);

\draw[fill=red!10!white, opacity = 0.4] (p1) -- (p4) -- ($(p4) + (0,-3*\h)$) -- ($(p1) + (0,-3*\h)$) -- (p1);

\begin{scope}[shift={(0.0,-3*\h)}]
\begin{scope}[shift = {($(0,-1.5*\h)$)}, every node/.append style={xslant=1},xslant=1, xscale= 1, yscale = 0.3]
\coordinate (p1) at (-3.0,-3.0);
\coordinate (p2) at (3.0,-3.0);
\coordinate (p3) at (3.0, 3.0);
\coordinate (p4) at (-3.0,3.0);
\draw[fill=yellow!50!white,opacity=0.4] (p1)--(p2)--(p3)--(p4)--(p1);
\end{scope}

\draw[fill=yellow!50!white, opacity = 0.4] (p1) --++ (0,3*\h) --++ (6.0,0) -- (p2) -- (p1);

\begin{scope}[shift = {($(0,1.5*\h)$)}, every node/.append style={xslant=1},xslant=1, xscale= 1, yscale = 0.3]
\coordinate (p1) at (-3.0,-3.0);
\coordinate (p2) at (3.0,-3.0);
\coordinate (p3) at (3.0, 3.0);
\coordinate (p4) at (-3.0,3.0);
\draw[fill=yellow!40!white,opacity=0.4] (p1)--(p2)--(p3)--(p4)--(p1);
\end{scope}

\draw[fill=yellow!50!white, opacity = 0.4] (p2) -- (p3) -- ($(p3) + (0,-3*\h)$) -- ($(p2) + (0,-3*\h)$) -- (p2);

\draw[fill=yellow!10!white, opacity = 0.4] (p1) -- (p4) -- ($(p4) + (0,-3*\h)$) -- ($(p1) + (0,-3*\h)$) -- (p1);

\end{scope}

\begin{scope}[shift={(6.0,0.0)}]
\begin{scope}[shift = {($(0,-1.5*\h)$)}, every node/.append style={xslant=1},xslant=1, xscale= 1, yscale = 0.3]
\coordinate (p1) at (-3.0,-3.0);
\coordinate (p2) at (3.0,-3.0);
\coordinate (p3) at (3.0, 3.0);
\coordinate (p4) at (-3.0,3.0);
\draw[fill=red!40!white,opacity=0.4] (p1)--(p2)--(p3)--(p4)--(p1);
\end{scope}

\draw[fill=red!40!white, opacity = 0.4] (p1) --++ (0,3*\h) --++ (6.0,0) -- (p2) -- (p1);

\begin{scope}[shift = {($(0,1.5*\h)$)}, every node/.append style={xslant=1},xslant=1, xscale= 1, yscale = 0.3]
\coordinate (p1) at (-3.0,-3.0);
\coordinate (p2) at (3.0,-3.0);
\coordinate (p3) at (3.0, 3.0);
\coordinate (p4) at (-3.0,3.0);
\draw[fill=red!40!white,opacity=0.4] (p1)--(p2)--(p3)--(p4)--(p1);
\end{scope}

\draw[fill=red!40!white, opacity = 0.4] (p2) -- (p3) -- ($(p3) + (0,-3*\h)$) -- ($(p2) + (0,-3*\h)$) -- (p2);

\draw[fill=red!10!white, opacity = 0.4] (p1) -- (p4) -- ($(p4) + (0,-3*\h)$) -- ($(p1) + (0,-3*\h)$) -- (p1);

\begin{scope}[shift={(0.0,-3*\h)}]
\begin{scope}[shift = {($(0,-1.5*\h)$)}, every node/.append style={xslant=1},xslant=1, xscale= 1, yscale = 0.3]
\coordinate (p1) at (-3.0,-3.0);
\coordinate (p2) at (3.0,-3.0);
\coordinate (p3) at (3.0, 3.0);
\coordinate (p4) at (-3.0,3.0);
\draw[fill=yellow!50!white,opacity=0.4] (p1)--(p2)--(p3)--(p4)--(p1);
\end{scope}

\draw[fill=yellow!50!white, opacity = 0.4] (p1) --++ (0,3*\h) --++ (6.0,0) -- (p2) -- (p1);

\begin{scope}[shift = {($(0,1.5*\h)$)}, every node/.append style={xslant=1},xslant=1, xscale= 1, yscale = 0.3]
\coordinate (p1) at (-3.0,-3.0);
\coordinate (p2) at (3.0,-3.0);
\coordinate (p3) at (3.0, 3.0);
\coordinate (p4) at (-3.0,3.0);
\draw[fill=yellow!50!white,opacity=0.4] (p1)--(p2)--(p3)--(p4)--(p1);
\end{scope}

\draw[fill=yellow!50!white, opacity = 0.4] (p2) -- (p3) -- ($(p3) + (0,-3*\h)$) -- ($(p2) + (0,-3*\h)$) -- (p2);

\draw[fill=yellow!10!white, opacity = 0.4] (p1) -- (p4) -- ($(p4) + (0,-3*\h)$) -- ($(p1) + (0,-3*\h)$) -- (p1);

\begin{scope}[shift = {($(-6.0,1.5*\h)$)}, every node/.append style={xslant=1},xslant=1, xscale= 1, yscale = 0.3]
\coordinate (p1) at (-1.0,-1.0);
\coordinate (p2) at (2*3.5,-1.0);
\coordinate (p3) at (2*3.5, 1.0);
\coordinate (p4) at (-1.0,1.0);
\draw[fill=black!70!white,opacity=0.4] (p1)--(p2)--(p3)--(p4)--(p1);
\end{scope}

\end{scope}

\end{scope}

\draw node at (1.0,2.6) {\scriptsize Unit cell $(1)$};
\draw node at (7.0,2.6) {\scriptsize Unit cell $(2)$};

\end{tikzpicture}
}}
\hfill 
\subfloat[\label{fig:sample_cell_equivalent_2} Equivalent problem]{\resizebox{0.45\columnwidth} {!}  {\begin{tikzpicture}[xscale= 0.25, yscale =0.4]
\pgfmathsetmacro{\h}{2*0.7}

\pgfmathsetmacro{\UnitcellWidth}{3.75}
\pgfmathsetmacro{\UnitcellLength}{3.75}
\pgfmathsetmacro{\UnitcellHeight}{3.75}
\pgfmathsetmacro{\numarrow}{3}
\pgfmathsetmacro{\numarrowa}{4}


\begin{scope}[shift = {($(0,-1.5*\h)$)}, every node/.append style={xslant=1},xslant=1, xscale= 1, yscale = 0.3]
\coordinate (p1) at (-3.0,-3.0);
\coordinate (p2) at (3.0,-3.0);
\coordinate (p3) at (3.0, 3.0);
\coordinate (p4) at (-3.0,3.0);
\draw[fill=red!40!white,opacity=0.4] (p1)--(p2)--(p3)--(p4)--(p1);
\end{scope}

\draw[fill=red!40!white, opacity = 0.4] (p1) --++ (0,3*\h) --++ (6.0,0) -- (p2) -- (p1);

\begin{scope}[shift = {($(0,1.5*\h)$)}, every node/.append style={xslant=1},xslant=1, xscale= 1, yscale = 0.3]
\coordinate (p1) at (-3.0,-3.0);
\coordinate (p2) at (3.0,-3.0);
\coordinate (p3) at (3.0, 3.0);
\coordinate (p4) at (-3.0,3.0);
\draw[fill=red!40!white,opacity=0.4] (p1)--(p2)--(p3)--(p4)--(p1);

\end{scope}

\draw[fill=red!40!white, opacity = 0.4] (p2) -- (p3) -- ($(p3) + (0,-3*\h)$) -- ($(p2) + (0,-3*\h)$) -- (p2);

\draw[fill=red!10!white, opacity = 0.4] (p1) -- (p4) -- ($(p4) + (0,-3*\h)$) -- ($(p1) + (0,-3*\h)$) -- (p1);

\begin{scope}[shift={(-1.5,-1.0)}]
\foreach \x in {0,...,\numarrow}
{
\foreach \y in {0,...,\numarrow}
{
	\coordinate (t1) at (-\UnitcellWidth*0.4 + \x*1.1*\UnitcellWidth/\numarrow, -\UnitcellLength*0.4 + \y*0.7*\UnitcellLength/\numarrow);
	\coordinate (t2) at (-\UnitcellWidth*0.4 + \x*1.1*\UnitcellWidth/\numarrow, -\UnitcellLength*0.4 + \y*0.7*\UnitcellLength/\numarrow + 0.18*\UnitcellLength);
\draw [blue, ->, thick] (t1) -- (t2);
}
}
\begin{scope}[shift={(0.55,0.0)}]
\foreach \x in {0,...,\numarrow}
{
\foreach \y in {0,...,\numarrow}
{
	\coordinate (t1) at (-\UnitcellWidth*0.4 + \x*1.1*\UnitcellWidth/\numarrow, -\UnitcellLength*0.4 + \y*0.7*\UnitcellLength/\numarrow);
	\coordinate (t2) at (-\UnitcellWidth*0.4 + \x*1.1*\UnitcellWidth/\numarrow, -\UnitcellLength*0.4 + \y*0.7*\UnitcellLength/\numarrow + 0.18*\UnitcellLength);
\draw [red, ->, thick] (t1) -- (t2);
}
}
\end{scope}

\end{scope}

\begin{scope}[shift = {($(-1.3,0.8*\h)$)}, every node/.append style={xslant=1},xslant=1, xscale= 1, yscale = 0.3]
\begin{scope}[shift={(-0.4,2.2)}]
\foreach \x in {0,...,\numarrow}
{
\foreach \y in {0,...,\numarrow}
{
	\coordinate (t1) at (-\UnitcellWidth*0.4 + \x*1.1*\UnitcellWidth/\numarrow, -\UnitcellLength*0.4 + \y*1.1*\UnitcellLength/\numarrow);
	\coordinate (t2) at (-\UnitcellWidth*0.4 + \x*1.1*\UnitcellWidth/\numarrow, -\UnitcellLength*0.4 + \y*1.1*\UnitcellLength/\numarrow + 0.28*\UnitcellLength);
\draw [blue, ->, thick] (t1) -- (t2);
}
}
\end{scope}

\begin{scope}[shift={(0.2,2.2)}]
\foreach \x in {0,...,\numarrow}
{
\foreach \y in {0,...,\numarrow}
{
	\coordinate (t1) at (-\UnitcellWidth*0.4 + \x*1.1*\UnitcellWidth/\numarrow, -\UnitcellLength*0.4 + \y*1.1*\UnitcellLength/\numarrow);
	\coordinate (t2) at (-\UnitcellWidth*0.4 + \x*1.1*\UnitcellWidth/\numarrow, -\UnitcellLength*0.4 + \y*1.1*\UnitcellLength/\numarrow + 0.28*\UnitcellLength);
\draw [red, ->, thick] (t1) -- (t2);
}
}
\end{scope}

\end{scope}

\begin{scope}[shift={(6.0,0.0)}]
\begin{scope}[shift = {($(0,-1.5*\h)$)}, every node/.append style={xslant=1},xslant=1, xscale= 1, yscale = 0.3]
\coordinate (p1) at (-3.0,-3.0);
\coordinate (p2) at (3.0,-3.0);
\coordinate (p3) at (3.0, 3.0);
\coordinate (p4) at (-3.0,3.0);
\draw[fill=red!40!white,opacity=0.4] (p1)--(p2)--(p3)--(p4)--(p1);
\end{scope}

\draw[fill=red!40!white, opacity = 0.4] (p1) --++ (0,3*\h) --++ (6.0,0) -- (p2) -- (p1);

\begin{scope}[shift = {($(0,1.5*\h)$)}, every node/.append style={xslant=1},xslant=1, xscale= 1, yscale = 0.3]
\coordinate (p1) at (-3.0,-3.0);
\coordinate (p2) at (3.0,-3.0);
\coordinate (p3) at (3.0, 3.0);
\coordinate (p4) at (-3.0,3.0);
\draw[fill=red!40!white,opacity=0.4] (p1)--(p2)--(p3)--(p4)--(p1);
\end{scope}

\draw[fill=red!40!white, opacity = 0.4] (p2) -- (p3) -- ($(p3) + (0,-3*\h)$) -- ($(p2) + (0,-3*\h)$) -- (p2);

\draw[fill=red!10!white, opacity = 0.4] (p1) -- (p4) -- ($(p4) + (0,-3*\h)$) -- ($(p1) + (0,-3*\h)$) -- (p1);

\begin{scope}[shift={(-1.5,-1.0)}]
\foreach \x in {0,...,\numarrow}
{
\foreach \y in {0,...,\numarrow}
{
	\coordinate (t1) at (-\UnitcellWidth*0.4 + \x*1.1*\UnitcellWidth/\numarrow, -\UnitcellLength*0.4 + \y*0.7*\UnitcellLength/\numarrow);
	\coordinate (t2) at (-\UnitcellWidth*0.4 + \x*1.1*\UnitcellWidth/\numarrow, -\UnitcellLength*0.4 + \y*0.7*\UnitcellLength/\numarrow + 0.18*\UnitcellLength);
\draw [blue, ->, thick] (t1) -- (t2);
}
}

\begin{scope}[shift={(0.55,0.0)}]
\foreach \x in {0,...,\numarrow}
{
\foreach \y in {0,...,\numarrow}
{
	\coordinate (t1) at (-\UnitcellWidth*0.4 + \x*1.1*\UnitcellWidth/\numarrow, -\UnitcellLength*0.4 + \y*0.7*\UnitcellLength/\numarrow);
	\coordinate (t2) at (-\UnitcellWidth*0.4 + \x*1.1*\UnitcellWidth/\numarrow, -\UnitcellLength*0.4 + \y*0.7*\UnitcellLength/\numarrow + 0.18*\UnitcellLength);
\draw [red, ->, thick] (t1) -- (t2);
}
}
\end{scope}

\end{scope}

\begin{scope}[shift = {($(-1.3,0.8*\h)$)}, every node/.append style={xslant=1},xslant=1, xscale= 1, yscale = 0.3]
\begin{scope}[shift={(-0.4,2.2)}]
\foreach \x in {0,...,\numarrow}
{
\foreach \y in {0,...,\numarrow}
{
	\coordinate (t1) at (-\UnitcellWidth*0.4 + \x*1.1*\UnitcellWidth/\numarrow, -\UnitcellLength*0.4 + \y*1.1*\UnitcellLength/\numarrow);
	\coordinate (t2) at (-\UnitcellWidth*0.4 + \x*1.1*\UnitcellWidth/\numarrow, -\UnitcellLength*0.4 + \y*1.1*\UnitcellLength/\numarrow + 0.28*\UnitcellLength);
\draw [blue, ->, thick] (t1) -- (t2);
}
}
\end{scope}

\begin{scope}[shift={(0.2,2.2)}]
\foreach \x in {0,...,\numarrow}
{
\foreach \y in {0,...,\numarrow}
{
	\coordinate (t1) at (-\UnitcellWidth*0.4 + \x*1.1*\UnitcellWidth/\numarrow, -\UnitcellLength*0.4 + \y*1.1*\UnitcellLength/\numarrow);
	\coordinate (t2) at (-\UnitcellWidth*0.4 + \x*1.1*\UnitcellWidth/\numarrow, -\UnitcellLength*0.4 + \y*1.1*\UnitcellLength/\numarrow + 0.28*\UnitcellLength);
\draw [red, ->, thick] (t1) -- (t2);
}
}
\end{scope}

\end{scope}

\begin{scope}[shift={(4.4,-0.4)}]
\foreach \x in {0,...,0}
{
\foreach \y in {0,...,\numarrow}
{
	\coordinate (t1) at (-\UnitcellWidth*0.4 + \x*0.8*\UnitcellWidth/\numarrow, -\UnitcellLength*0.4 + \y*0.7*\UnitcellLength/\numarrow);
	\coordinate (t2) at (-\UnitcellWidth*0.4 + \x*0.8*\UnitcellWidth/\numarrow, -\UnitcellLength*0.4 + \y*0.7*\UnitcellLength/\numarrow + 0.18*\UnitcellLength);
\draw [blue, ->, thick] (t1) -- (t2);
}
}
\end{scope}

\begin{scope}[shift={(4.85,-0.1)}]
\foreach \x in {0,...,0}
{
\foreach \y in {0,...,\numarrow}
{
	\coordinate (t1) at (-\UnitcellWidth*0.4 + \x*0.8*\UnitcellWidth/\numarrow, -\UnitcellLength*0.4 + \y*0.7*\UnitcellLength/\numarrow);
	\coordinate (t2) at (-\UnitcellWidth*0.4 + \x*0.8*\UnitcellWidth/\numarrow, -\UnitcellLength*0.4 + \y*0.7*\UnitcellLength/\numarrow + 0.18*\UnitcellLength);
\draw [red, ->, thick] (t1) -- (t2);
}
}
\end{scope}

\end{scope}

\draw node at (1.0,3.6) {\scriptsize $\mathcal S_{\text{eq}}^{(1)}$};
\draw node at (7.0,3.6) {\scriptsize $\mathcal S_{\text{eq}}^{(2)}$};

\end{tikzpicture}
}}
\hfill \null
\caption{(a) An example of an array composed of two unit cells featuring connected  PEC traces. Each unit cell has two dielectric regions (shown with yellow and red colors) and is enclosed by a fictitious surface. The PEC trace is cut into two pieces by fictitious surfaces $\mathcal S_{\text{eq}}^{(1)}$ and $\mathcal S_{\text{eq}}^{(2)}$. (b) Unit cells are modeled by equivalent surface electric and magnetic current densities (shown with blue and red arrows, respectively) introduced on $\mathcal S_{\text{eq}}^{(1)}$ and $\mathcal S_{\text{eq}}^{(2)}$.}
\label{fig:equivalent_cell}
\end{figure}
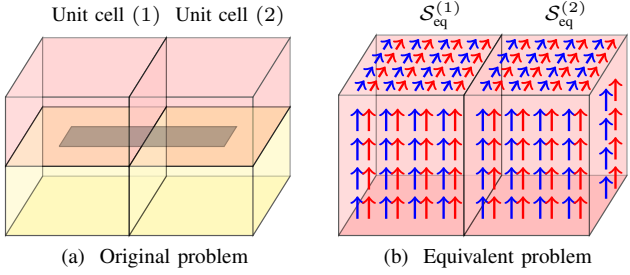

\section{Proposed Method}
\label{proposed method}
A main limitation of the previous work~\cite{Utkarsh_macromodel_2020} is that the conductor traces are not allowed to traverse the equivalent surfaces. Therefore, the efficient simulation of electromagnetic surfaces that are composed of slot shapes and are employed to tailor the electric and magnetic fields, especially in transmission mode, remains unaddressed. Hence, accurate modeling of the continuous current of the conductor traces cut by the equivalent surfaces in macromodeling algorithms constitutes the core contribution of this paper in relation to~\cite{Utkarsh_macromodel_2020}. In this section, we explain how to handle this problem efficiently. 

\subsection{Surface Integral Equation}
In order to simplify the explanation of the proposed technique for analysis of unit cells with connected conductors, we consider a time-harmonic scattering problem from an EM surface consisting of two dielectric layers situated in free space, such as the one shown in Fig.~\ref{fig:equivalent_cell}\subref{fig:sample_cell_equivalent_1}. Let us assume a rectangular perfect electric conductor (PEC) trace between two adjacent unit cells. 
Similar to the previous work, we first enclose the $m$-th unit cell with a fictitious surface $\mathcal S_\text{eq}^{(m)}$, as shown in Fig.~\ref{fig:equivalent_cell} ($m = 1, 2$ here).
The fictitious surfaces intersect with the PEC trace, and form the boundaries of each macromodel.

\begin{figure}
\centering
\begin{tikzpicture}[scale=.8]

\filldraw[fill=red!40!white, draw=black] (0,0) rectangle (4,1.8);
\filldraw[fill=yellow!70!white, draw=black] (0,-1.8) rectangle (4,0);
\draw [line width=0.9mm] (1,0) -- (4,0);
\draw node at (2.0,2.2) {$\mathcal S_{\text{eq}}^{(1)}$};
\draw node at (2.0,1.0) {$\mathcal{V}_2^{(1)}$};
\draw node at (2.0,-1.0) {$\mathcal{V}_1^{(1)}$};
\draw[red,thick] (0,-1.8) rectangle (4.0,1.8);

\begin{scope}[shift={(+4.0,0.0)}]

\filldraw[fill=red!40!white, draw=black] (0,0) rectangle (4,1.8);
\filldraw[fill=yellow!70!white, draw=black] (0,-1.8) rectangle (4,0);
\draw [line width=0.9mm] (1-1,0) -- (4-1,0);
\draw node at (2.0,2.2) {$\mathcal S_{\text{eq}}^{(2)}$};
\draw node at (2.0,1.0) {$\mathcal{V}_2^{(2)}$};
\draw node at (2.0,-1.0) {$\mathcal{V}_1^{(2)}$};
\draw[red,thick] (0,-1.8) rectangle (4.0,1.8);

\end{scope}

\end{tikzpicture}
\caption{Side view of two adjacent unit cells considered in Fig.~\ref{fig:equivalent_cell} with dielectric regions $\mathcal V_1^{(m)}$ and $\mathcal V_2^{(m)}$, $m=1,2$. The left and right unit cells are enclosed by fictitious surfaces $\mathcal S_{\text{eq}}^{(1)}$ and $\mathcal S_{\text{eq}}^{(2)}$, respectively,  which are drawn in red. }
\label{fig:unit cells side view}
\end{figure}
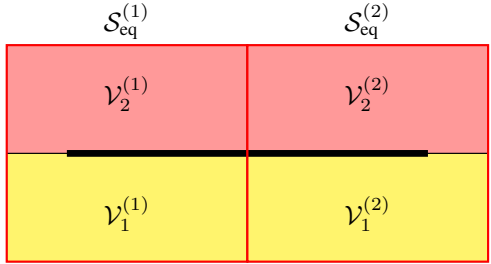
The cross-section of the two unit cells are depicted in Fig.~\ref{fig:unit cells side view} with the PEC trace in between, and the dielectric regions denoted by $\mathcal V_v^{(m)}$, with region index $v$ ($v=1, 2$). The surface enclosing the $v$-th region is $\mathcal S_{v}^{(m)}$. Next, $\mathcal S_{v}^{(m)}$ is discretized with triangular elements. For better visualization of the discretization step, we separate the surfaces that need to be discretized for each unit cell in Fig.~\ref{fig:unit cell discretization scheme}. Similar to the Section~\ref{Surface Integral Equation}, we apply the EFIE and the MFIE to each region to obtain the system of matrix~\eqref{SIE Matrix equation-1} for a given unit cell.  
  

\begin{figure}
\centering
\begin{tikzpicture}[scale=0.9]

\filldraw[fill=red!40!white, draw=black] (0,0+.3) rectangle (4,1.8+.3);
\draw[black,thick,dashed] (0+.1,0+.3+.1) rectangle (4-.1,1.8+.3-.1);
\filldraw[fill=yellow!70!white, draw=black] (0,-1.8-.3) rectangle (4,0-.3);
\draw[black,thick,dashed] (0+.1,-1.8-.3+.1) rectangle (4-.1,0-.3-.1);
\draw [line width=0.9mm] (1,0) -- (4,0);
\draw[red,thick,dashed] (-0.5,-2.6) rectangle (4.5,2.6);
\draw node at (2.0,1.0+.2) {$\mathcal{V}_2^{(m)}$};
\draw node at (2.0,-1.0-.2) {$\mathcal{V}_1^{(m)}$};
\draw node at (1.2,1.0+.65) {$\mathcal S_2^{(m)}$};
\draw node at (1.2,-1.0-.7) {$\mathcal S_1^{(m)}$};

\draw node at (-.98,2.2) {$\mathcal S_{\text{eq}}^{(m)}$};





\end{tikzpicture}
\caption{Discretization of the surfaces of a unit cell featuring connected PEC traces. Regions $\mathcal V_1^{(m)}$ and $\mathcal V_2^{(m)}$ are enclosed by surfaces $\mathcal S_1^{(m)}$ and $\mathcal S_2^{(m)}$, respectively. The PEC trace and the regions are separated in this figure only for visualization clarity.}
\label{fig:unit cell discretization scheme}
\end{figure}
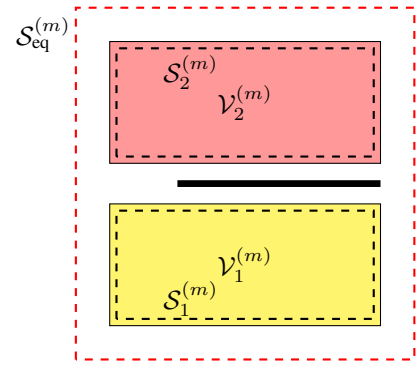

\begin{figure}
\centering
\begin{tikzpicture}[scale=1]

\filldraw[fill=red!40!white, draw=black] (0,0) rectangle (4,1.8);
\filldraw[fill=yellow!70!white, draw=black] (0,-1.8) rectangle (4,0);
\draw [line width=0.9mm] (1,0) -- (4,0);
\draw [line width=0.9mm,dashed] (4,0) -- (5,0);
\draw[red,thick] (0,-1.8) rectangle (4.0,1.8);
\draw node at (2.0,1.0) {$\mathcal{V}_2^{(m)}$};
\draw node at (2.0,-1.0) {$\mathcal{V}_1^{(m)}$};
\draw node at (2.0,2.1) {$\mathcal S_{\text{eq}}^{(m)}$};

\draw (4.3,0.6) rectangle (4.3+.4,0.6+.4) node at (4.3+.2,0.6+.2) {$1$};
\draw node at (4.4+.2,0.2+.1) {$j_{\text{eq},n}$};
\draw node at (4.-.5,0.25+.1) {$j_{2,\hat n}$};
\draw node at (4.4+.25,-0.4) {$j_{\text{eq},n'}$};
\draw node at (4.-.5,-0.35) {$j_{1,\Tilde{n}}$};
\draw [red][line width=0.45mm,-latex] (4.15,0.6) -- (4.15,0+0.1);
\draw [red][line width=0.45mm,-latex] (4.15,0-0.1) -- (4.15,-0.6);

\draw [blue][line width=0.45mm] (3.4+.1,0.1) -- (3.8+.1,0.1);
\draw [blue][line width=0.45mm,-latex] (3.78+.1,0.1) -- (3.78+.1,0.1+.5);

\draw [blue][line width=0.45mm,-latex]  (3.9,-.15) -- (3.45,-.15);
\draw [blue][line width=0.45mm] (3.88,-0.15) -- (3.88,-0.56);

\draw (-.7,0.3) rectangle (-.7+.4,0.3+.4) node at (-.7+.2,0.3+.2) {$2$} ;
\draw [blue][line width=0.45mm,-latex] (-0.15,-0.4) -- (-0.15,0.48);

\draw [blue][line width=0.45mm,-latex] (0.08,0.1) -- (0.5,0.1);
\draw [blue][line width=0.45mm] (0.1,0.1) -- (0.1,0.1+.35);

\draw [blue][line width=0.45mm] (0.1,-.08) -- (0.45,-.08);
\draw [blue][line width=0.45mm,-latex] (0.1,-0.06) -- (0.1,-0.48);

\draw (1.5,+.2) rectangle (1.5+.4,+.2+.4) node at (1.5+.2,+.2+.2) {$3$} ;
\draw [blue][line width=0.45mm,-latex] (0.75,0.15) -- (1.5,0.15);
\draw [blue][line width=0.45mm,-latex] (0.75,-0.15) -- (1.5,-0.15);

\draw (2.6,-1.7) rectangle (2.6+.4,-1.7+.4) node at (2.6+.2,-1.7+.2) {$4$} ;
\draw [blue][line width=0.45mm,-latex] (1.75,-1.65) -- (2.55,-1.65);
\draw [blue][line width=0.45mm,-latex] (1.75,-1.95) -- (2.55,-1.95);

\end{tikzpicture}
\caption{Special junctions for enforcement of the boundary conditions are labeled $\boxed{1}, \dots, \boxed{4}$. Half RWG basis functions are shown in red and full RWG basis functions are shown in blue. The dashed line shows where the PEC trace traverses the equivalent surface $\mathcal S_{\text{eq}}^{(m)}$.}
\label{fig:Boundary_Condition_Interior}
\end{figure}
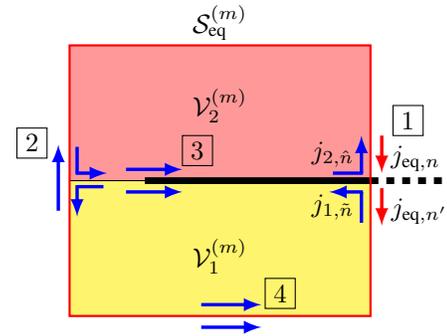

\subsection{Enforcement of the interior Boundary Conditions and Interelement Coupling}
\label{Enforcement of Boundary Conditions}
Accurate modeling of the continuous current flowing in the PEC trace cut by the two adjacent equivalent surfaces can be obtained through the enforcement of the boundary conditions. In this section, we will discuss the relevant boundary conditions. Here, a junction is defined to be an intersection of three or more regions, such as junctions that are labeled with $\boxed{1}$ and $\boxed{2}$ in Fig~\ref{fig:Boundary_Condition_Interior}. Two scenarios can arise when a PEC trace intersects with the equivalent surface $\mathcal S_{\text{eq}}$, depending on whether the trace extends into the next unit cell or terminates right on $\mathcal S_{\text{eq}}$.     

\begin{enumerate}
\label{Enforcement of Boundary Condition1}
\item{\bf PEC trace intersecting $\mathcal S_{\text{eq}}$ and extending into the next unit cell:} Consider the sample junction on $\mathcal S_{\text{eq}}^{(m)}$ that is labeled with $\boxed{1}$ where the PEC trace, after intersecting $\mathcal S_{\text{eq}}$, continues into the next unit cell on the right. Due to the discontinuity of the tangential magnetic field on the two sides of the PEC trace between regions $\mathcal{V}_1^{(m)}$ and $\mathcal{V}_2^{(m)}$, we need to introduce two independent electric current densities, as shown with blue arrows in Fig.~\ref{fig:Boundary_Condition_Interior}. These two currents are defined on the interface of the regions and expanded with full RWG basis functions with corresponding coefficients $j_{1,\Tilde{n}}$ and $j_{2,\hat n}$. Likewise, we need to define two independent electric currents on $\mathcal S_{\text{eq}}^{(m)}$. These two currents are expanded with half RWG basis functions shown with red arrows in Fig.~\ref{fig:Boundary_Condition_Interior}, with corresponding coefficients $j_{\text{eq},n}$ and $j_{\text{eq},n'}$. In addition, the tangential magnetic field is continuous across the equivalent surface $\mathcal S_{\text{eq}}$. Hence, considering the orientation of the basis functions, we have $j_{2,\hat n}=j_{\text{eq},n}$ and $j_{1,\Tilde{n}} = j_{\text{eq},n'}$. We collect $j_{2,\hat n}$ and $j_{1,\Tilde{n}}$ into the vector $\widetilde{\textbf X}_{\text{eq}}^{(m)}$. To enforce these boundary conditions using the sparse matrix $\textbf U^{(m)}$ in~\eqref{UX}, the entries $(q^{j_{2,\hat n}},\Tilde{q}^{j_{\text{eq},n}})$ and $(q^{j_{1,\Tilde n}},\Tilde{q}^{j_{\text{eq},n'}})$ of $\textbf U^{(m)}$ need to be set to $1$, where $q^{\alpha}$ and $\Tilde{q}^{\beta}$ are, respectively, the indices of the entries associated with the coefficient $\alpha$ in $\textbf X^{(m)}$ and $\beta$ in $\widetilde{\textbf X}^{(m)}$. There is no need to discretize the magnetic current $\vec M_{\text{eq}}$ on the edge of the PEC trace at the intersection of $\mathcal S_{\text{eq}}$ since the electric field tangential to the PEC trace is zero.

\item{\bf PEC trace terminating on $\mathcal S_{\text{eq}}$:} If the PEC trace in~Fig.~\ref{fig:Boundary_Condition_Interior} is not connected to another conductive element from the right side, the current density normal to $\mathcal S_{\text{eq}}$ at the junction labeled with $\boxed{1}$ is zero. Therefore, the tangential magnetic field is continuous on the two sides of the PEC interface between regions $\mathcal{V}_1^{(m)}$ and $\mathcal{V}_2^{(m)}$, which leads to $j_{1,\Tilde{n}}=j_{2,\hat n}$. The continuity of the tangential magnetic field is enforced similarly on $\mathcal S_{\text{eq}}^{(m)}$ which leads to $j_{\text{eq},n}=j_{\text{eq},n'}$. In this case, the two half RWG basis functions become equivalent to a full RWG basis function. Moreover, due to the continuity of the tangential magnetic field across $\mathcal S_{\text{eq}}$, we have $j_{2,\hat n}=j_{1,\Tilde{n}}=j_{\text{eq},n}=j_{\text{eq},n'}$. Therefore, only one electric current unknown is needed for this edge, with coefficient $j_{\text{eq},n^{''}}$, which we collect into the vector $\widetilde{\textbf X}_{\text{eq}}^{(m)}$. To enforce these boundary conditions in \eqref{UX}, we can set the entries $(q^{j_{2,\hat{n}}},\Tilde{q}^{j_{\text{eq},n''}})$, $(q^{j_{1,\Tilde n}},\Tilde{q}^{j_{\text{eq},n''}})$, $(q^{j_{\text{eq},n}},\Tilde{q}^{j_{\text{eq},n''}})$, and $(q^{j_{\text{eq},n'}},\Tilde{q}^{j_{\text{eq},n''}})$ of $\textbf U^{(m)}$ to be $1$. 
\end{enumerate}  

The other edges and junctions that are labeled with $\boxed{2}$, $\boxed{3}$, and $\boxed{4}$ in Fig.~\ref{fig:Boundary_Condition_Interior} are treated according to~\cite{Utkarsh_macromodel_2020}. Next, we generate the macromodeling matrix of the $m$-th unit cell $\widetilde{\textbf Z}_{\text{eq},\text{eq}}^{(m)}$ using \eqref{Shur_complement}. The above macromodeling matrix generation step can be generalized for a unit cell consisting of an arbitrary number of dielectric layers with arbitrary shaped PEC scatterers traversing the macromodel boundaries through multiple layers in a straightforward manner. 

\begin{figure}
\centering
\begin{tikzpicture}[scale=1]

\draw[red,thick] (0,-1.8) rectangle (3.0,1.8);
\draw node at (1.5,2.1) {$\mathcal S_{\text{eq}}^{(m)}$};
\draw (3.1,-1.1) rectangle (3.1+.4,-1.1+.4) node at (3.1+.2,-1.1+.2) {$1$};
\draw node at (2.35,0.2+.2) {$j_{\text{eq},n}^{(m)}$};
\draw node at (2.35,-0.5) {$j_{\text{eq},\hat{n}}^{(m)}$};
\draw [line width=0.9mm,dashed] (2.0,0) -- (3,0);
\draw [red][line width=0.45mm,-latex] (3.15-.3,0.6) -- (3.15-.3,0+0.1);
\draw [red][line width=0.45mm,-latex] (3.15-.3,0-0.1) -- (3.15-.3,-0.6);

\begin{scope}[shift={(+3.6,0.0)}]
\draw[red,thick] (0,-1.8) rectangle (3.0,1.8);
\draw node at (1.5,2.1) {$\mathcal S_{\text{eq}}^{(m')}$};
\draw node at (0.45+.2,0.2+.2) {$j_{\text{eq},n'}^{(m')}$};
\draw node at (0.45+.25,-0.5) {$j_{\text{eq},\Tilde{n}}^{(m')}$};
\draw [line width=0.9mm,dashed] (0.0,0) -- (1.0,0);

\draw [red][line width=0.45mm,-latex] (0.15,0+0.1) -- (0.15,0.6);
\draw [red][line width=0.45mm,-latex] (0.15,-0.6) -- (0.15,0-0.1);
\end{scope}

\draw[black,thick,dashed] (3.0,1.8) -- (3.0,1.8+.6);
\draw[black,thick,dashed] (3.6,1.8) -- (3.6,1.8+.6);
\draw[black,thick,<->] (3.0,1.8+.6) -- (3.6,1.8+.6);
\draw node at (3.3,1.8+.9) {$\Delta \rightarrow 0$};

\begin{scope}[shift={(-1.0,1.8)}]
\draw [blue][line width=0.45mm,-latex]  (3.9,-.15) -- (3.45,-.15);
\draw [blue][line width=0.45mm] (3.88,-0.15) -- (3.88,-0.56);
\end{scope}

\draw (3.1,1.3) rectangle (3.1+.4,1.3+.4) node at (3.1+.2,1.3+.2) {$2$};
\begin{scope}[shift={(3.65,1.73)}]
\draw [blue][line width=0.45mm] (0.1,-.08) -- (0.45,-.08);
\draw [blue][line width=0.45mm,-latex] (0.1,-0.06) -- (0.1,-0.53);
\end{scope}

\end{tikzpicture}
\caption{Special junctions for enforcement of the boundary conditions when two equivalent surfaces are connected. The space between the surfaces is only shown here for clear representation of the current directions and is actually zero. Half RWG basis functions are associated with the edges of the PEC traces traversing the equivalent surfaces and are shown in red. Full RWG basis functions are associated with regular edges and are shown in blue. Dashed lines show the positions where PEC traces traverse the equivalent surfaces.}
\label{fig:Boundary Condition Exterior}
\end{figure}
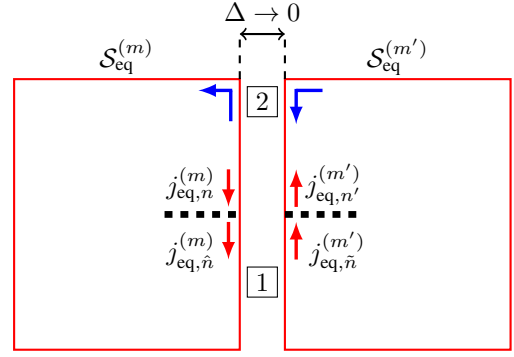

Next, we model the interelement coupling between the equivalent surfaces. An example of two adjacent equivalent surfaces $S^{(m)}_\text{eq}$ and $S^{(m')}_\text{eq}$ with PEC trace traversing the equivalent surfaces is depicted in Fig.~\ref{fig:Boundary Condition Exterior}. The interaction between the unit cells can be computed through the equivalent current densities on the macromodel boundaries using~\eqref{SLAE_exterior EFIE-MFIE-1}. The following boundary condition needs to be enforced between $S^{(m)}_\text{eq}$ and $S^{(m')}_\text{eq}$ through matrix $\textbf U_0$ in~\eqref{U0Y}: 

\begin{itemize}
\item{\bf Intersection of a PEC trace of adjacent unit cells:} We consider the edge of a PEC trace on the equivalent surfaces that is labeled with $\boxed{1}$ in Fig.~\ref{fig:Boundary Condition Exterior}. As we discussed earlier in this section, two independent electric currents need to be introduced on each equivalent surface to accurately model the discontinuity in the tangential magnetic field that arises when a current flows on the PEC trace. These two currents are expanded by half RWG basis functions with coefficients $j_{\text{eq},n}^{(m)}$ and $j_{\text{eq},\hat{n}}^{(m)}$ on $S^{(m)}_\text{eq}$ and coefficients $j_{\text{eq},n'}^{(m')}$ and $j_{\text{eq},\Tilde{n}}^{(m')}$ on $S^{(m')}_\text{eq}$. Since the tangential magnetic field remains continuous across the interface between the two equivalent surfaces, the electric current density coefficient $j_{\text{eq},n}^{(m)}$ on $S^{(m)}_\text{eq}$ is set to be equal to $j_{\text{eq},{n'}}^{(m')}$ on $S^{(m')}_\text{eq}$. Likewise, the electric current density coefficient $j_{\text{eq},\hat{n}}^{(m)}$ on $S^{(m)}_\text{eq}$ is set to be equal to $j_{\text{eq},\Tilde{n}}^{(m')}$ on $S^{(m')}_\text{eq}$. To enforce these boundary conditions by~\eqref{U0Y}, the entries $(q^{j_{\text{eq},n'}^{(m')}},\Tilde{q}^{j_{\text{eq},{n}}^{(m)}})$ and $(q^{j_{\text{eq},\Tilde{n}}^{(m')}},\Tilde{q}^{j_{\text{eq},\hat{n}}^{(m)}})$ of $\textbf U_0$ need to be set to $1$.


\end{itemize}
A detailed description of other edges and junctions that are labeled with $\boxed{2}$ in Fig.~\ref{fig:Boundary Condition Exterior} is provided in~\cite{Utkarsh_macromodel_2020}.

\section{FFT Acceleration}
The final system of equations in the proposed approach can now be expressed in a form similar to~\eqref{U0tZU0X_matrix-1} by applying the appropriate boundary conditions described above. To enable the simulation of large-scale surfaces, we study the extension of the previously developed FFT-based acceleration technique~\cite{Utkarsh_macromodel_2020} to the proposed approach. 

\begin{figure}[t]
\null \hfill
\subfloat[\label{fig:conductor_traverse_meshing1} ] {\resizebox{1.03\columnwidth} {!} {\begin{tikzpicture}[scale=1]

\draw[black] node at (1.0,2.15) {\scriptsize{$\mathcal S_{\text{eq}}^{(1)}$}};
\draw[black,thick] (0,0) rectangle (1.8,1.8);
\draw[red] (0.3,0.5) -- (0.3+0.3,0.5);
\draw[red] (0.3+0.3,0.5) -- (0.3+0.15,0.5/2);
\draw[red] (0.3+0.15,0.5/2) -- (0.3,0.5);

\draw[red] (0.3,0.5) -- (0.3+0.3,0.5);
\draw[red] (0.3+0.3,0.5) -- (0.3+0.15,0.5+0.25);
\draw[red] (0.3+0.15,0.5+0.25) -- (0.3,0.5);
\draw[red] node at (0.3+0.6,0.5) {\scriptsize{$1,a$}};

\draw[blue] (1.0,1.5) -- (1.0+0.3,1.5+0.15);
\draw[blue] (1.0+0.3,1.5+0.15) -- (1.0+0.3,1.5-0.15);
\draw[blue] (1.0+0.3,1.5-0.15) -- (1.0,1.5);

\draw[blue] (1.0+0.3,1.5+0.15) -- (1.0+0.3,1.5-0.15);
\draw[blue] (1.0+0.3,1.5-0.15) -- (1.6,1.5);
\draw[blue] (1.6,1.5) -- (1.0+0.3,1.5+0.15);
\draw[blue] node at (1.0+0.3,1.5-0.28) {\scriptsize{$2,b$}};

\begin{scope}[shift={(0.53,-0.65)}]
\draw[cyan] (1.0+0.3,1.5+0.15) -- (1.0+0.3,1.5-0.15);
\draw[cyan] (1.0+0.3,1.5-0.15) -- (1.6,1.5);
\draw[cyan] (1.6,1.5) -- (1.0+0.3,1.5+0.15);
\draw[cyan] node at (1.55,1.8) {\scriptsize{$3,c$}};
\draw [line width=0.9mm,dashed] (0.7,1.5) -- (1.3,1.5);
\end{scope}

\begin{scope}[shift={(-1.33,-0.65)}]
\draw[cyan] (1.0,1.5) -- (1.0+0.3,1.5+0.15);
\draw[cyan] (1.0+0.3,1.5+0.15) -- (1.0+0.3,1.5-0.15);
\draw[cyan] (1.0+0.3,1.5-0.15) -- (1.0,1.5);
\draw[cyan] node at (1.6,1.75) {\scriptsize{$4,d$}};
\draw [line width=0.9mm,dashed] (1.35,1.5) -- (1.35+0.6,1.5);
\end{scope}

\begin{scope}[shift={(+2.8,0.0)}]
\draw[black] node at (1.0,2.15) {\scriptsize{$\mathcal S_{\text{eq}}^{(2)}$}};
\draw[black,thick] (0,0) rectangle (1.8,1.8);
\draw[red] (0.3,0.5) -- (0.3+0.3,0.5);
\draw[red] (0.3+0.3,0.5) -- (0.3+0.15,0.5/2);
\draw[red] (0.3+0.15,0.5/2) -- (0.3,0.5);

\draw[red] (0.3,0.5) -- (0.3+0.3,0.5);
\draw[red] (0.3+0.3,0.5) -- (0.3+0.15,0.5+0.25);
\draw[red] (0.3+0.15,0.5+0.25) -- (0.3,0.5);
\draw[red] node at (0.3+0.6,0.5) {\scriptsize{$5,a$}};

\draw[blue] (1.0,1.5) -- (1.0+0.3,1.5+0.15);
\draw[blue] (1.0+0.3,1.5+0.15) -- (1.0+0.3,1.5-0.15);
\draw[blue] (1.0+0.3,1.5-0.15) -- (1.0,1.5);

\draw[blue] (1.0+0.3,1.5+0.15) -- (1.0+0.3,1.5-0.15);
\draw[blue] (1.0+0.3,1.5-0.15) -- (1.6,1.5);
\draw[blue] (1.6,1.5) -- (1.0+0.3,1.5+0.15);
\draw[blue] node at (1.0+0.3,1.5-0.28) {\scriptsize{$6,b$}};

\begin{scope}[shift={(0.53,-0.65)}]
\draw[cyan] (1.0+0.3,1.5+0.15) -- (1.0+0.3,1.5-0.15);
\draw[cyan] (1.0+0.3,1.5-0.15) -- (1.6,1.5);
\draw[cyan] (1.6,1.5) -- (1.0+0.3,1.5+0.15);
\draw[cyan] node at (1.55,1.8) {\scriptsize{$7,c$}};
\draw [line width=0.9mm,dashed] (0.7,1.5) -- (1.3,1.5);
\end{scope}

\begin{scope}[shift={(-1.33,-0.65)}]
\draw[cyan] (1.0,1.5) -- (1.0+0.3,1.5+0.15);
\draw[cyan] (1.0+0.3,1.5+0.15) -- (1.0+0.3,1.5-0.15);
\draw[cyan] (1.0+0.3,1.5-0.15) -- (1.0,1.5);
\draw[cyan] node at (1.6,1.75) {\scriptsize{$8,d$}};
\draw [line width=0.9mm,dashed] (1.35,1.5) -- (1.35+0.6,1.5);
\end{scope}

\end{scope}

\begin{scope}[shift={(+2.8*2,0.0)}]
\draw[black] node at (1.0,2.15) {\scriptsize{$\mathcal S_{\text{eq}}^{(3)}$}};
\draw[black,thick] (0,0) rectangle (1.8,1.8);
\draw[red] (0.3,0.5) -- (0.3+0.3,0.5);
\draw[red] (0.3+0.3,0.5) -- (0.3+0.15,0.5/2);
\draw[red] (0.3+0.15,0.5/2) -- (0.3,0.5);

\draw[red] (0.3,0.5) -- (0.3+0.3,0.5);
\draw[red] (0.3+0.3,0.5) -- (0.3+0.15,0.5+0.25);
\draw[red] (0.3+0.15,0.5+0.25) -- (0.3,0.5);
\draw[red] node at (0.3+0.6,0.5) {\scriptsize{$9,a$}};

\draw[blue] (1.0,1.5) -- (1.0+0.3,1.5+0.15);
\draw[blue] (1.0+0.3,1.5+0.15) -- (1.0+0.3,1.5-0.15);
\draw[blue] (1.0+0.3,1.5-0.15) -- (1.0,1.5);

\draw[blue] (1.0+0.3,1.5+0.15) -- (1.0+0.3,1.5-0.15);
\draw[blue] (1.0+0.3,1.5-0.15) -- (1.6,1.5);
\draw[blue] (1.6,1.5) -- (1.0+0.3,1.5+0.15);
\draw[blue] node at (1.0+0.3,1.5-0.28) {\scriptsize{$10,b$}};

\begin{scope}[shift={(0.53,-0.65)}]
\draw[cyan] (1.0+0.3,1.5+0.15) -- (1.0+0.3,1.5-0.15);
\draw[cyan] (1.0+0.3,1.5-0.15) -- (1.6,1.5);
\draw[cyan] (1.6,1.5) -- (1.0+0.3,1.5+0.15);
\draw[cyan] node at (1.6,1.8) {\scriptsize{$11,c$}};
\draw [line width=0.9mm,dashed] (0.7,1.5) -- (1.3,1.5);
\end{scope}

\begin{scope}[shift={(-1.33,-0.65)}]
\draw[cyan] (1.0,1.5) -- (1.0+0.3,1.5+0.15);
\draw[cyan] (1.0+0.3,1.5+0.15) -- (1.0+0.3,1.5-0.15);
\draw[cyan] (1.0+0.3,1.5-0.15) -- (1.0,1.5);
\draw[cyan] node at (1.65,1.75) {\scriptsize{$12,d$}};
\draw [line width=0.9mm,dashed] (1.35,1.5) -- (1.35+0.6,1.5);
\end{scope}

\end{scope}

\end{tikzpicture}}}
\hfill \null
\\
\null \hfill
\subfloat[\label{fig:conductor_traverse_meshing2} ] {\resizebox{0.92\columnwidth} {!} {\begin{tikzpicture}[scale=1]

\draw[black] node at (1.0,2.15) {\scriptsize{$\mathcal S_{\text{eq}}^{(1)}$}};
\draw[black,thick] (0,0) rectangle (1.8,1.8);
\draw[red] (0.3,0.5) -- (0.3+0.3,0.5);
\draw[red] (0.3+0.3,0.5) -- (0.3+0.15,0.5/2);
\draw[red] (0.3+0.15,0.5/2) -- (0.3,0.5);

\draw[red] (0.3,0.5) -- (0.3+0.3,0.5);
\draw[red] (0.3+0.3,0.5) -- (0.3+0.15,0.5+0.25);
\draw[red] (0.3+0.15,0.5+0.25) -- (0.3,0.5);
\draw[red] node at (0.3+0.6,0.5) {\scriptsize{$1,a$}};

\draw[blue] (1.0,1.5) -- (1.0+0.3,1.5+0.15);
\draw[blue] (1.0+0.3,1.5+0.15) -- (1.0+0.3,1.5-0.15);
\draw[blue] (1.0+0.3,1.5-0.15) -- (1.0,1.5);

\draw[blue] (1.0+0.3,1.5+0.15) -- (1.0+0.3,1.5-0.15);
\draw[blue] (1.0+0.3,1.5-0.15) -- (1.6,1.5);
\draw[blue] (1.6,1.5) -- (1.0+0.3,1.5+0.15);
\draw[blue] node at (1.0+0.3,1.5-0.25) {\scriptsize{$2,b$}};

\begin{scope}[shift={(0.53,-0.65)}]
\draw[cyan] (1.0+0.3,1.5+0.15) -- (1.0+0.3,1.5-0.15);
\draw[cyan] (1.0+0.3,1.5-0.15) -- (1.6,1.5);
\draw[cyan] (1.6,1.5) -- (1.0+0.3,1.5+0.15);
\draw[cyan] node at (1.55,1.8) {\scriptsize{$3,c$}};
\draw [line width=0.9mm,dashed] (0.7,1.5) -- (1.3,1.5);
\end{scope}

\begin{scope}[shift={(+2.8,0.0)}]
\draw[black] node at (1.0,2.15) {\scriptsize{$\mathcal S_{\text{eq}}^{(2)}$}};
\draw[black,thick] (0,0) rectangle (1.8,1.8);
\draw[red] (0.3,0.5) -- (0.3+0.3,0.5);
\draw[red] (0.3+0.3,0.5) -- (0.3+0.15,0.5/2);
\draw[red] (0.3+0.15,0.5/2) -- (0.3,0.5);

\draw[red] (0.3,0.5) -- (0.3+0.3,0.5);
\draw[red] (0.3+0.3,0.5) -- (0.3+0.15,0.5+0.25);
\draw[red] (0.3+0.15,0.5+0.25) -- (0.3,0.5);
\draw[red] node at (0.3+0.6,0.5) {\scriptsize{$4,a$}};

\draw[blue] (1.0,1.5) -- (1.0+0.3,1.5+0.15);
\draw[blue] (1.0+0.3,1.5+0.15) -- (1.0+0.3,1.5-0.15);
\draw[blue] (1.0+0.3,1.5-0.15) -- (1.0,1.5);

\draw[blue] (1.0+0.3,1.5+0.15) -- (1.0+0.3,1.5-0.15);
\draw[blue] (1.0+0.3,1.5-0.15) -- (1.6,1.5);
\draw[blue] (1.6,1.5) -- (1.0+0.3,1.5+0.15);
\draw[blue] node at (1.0+0.3,1.5-0.25) {\scriptsize{$5,b$}};

\begin{scope}[shift={(0.53,-0.65)}]
\draw[cyan] (1.0+0.3,1.5+0.15) -- (1.0+0.3,1.5-0.15);
\draw[cyan] (1.0+0.3,1.5-0.15) -- (1.6,1.5);
\draw[cyan] (1.6,1.5) -- (1.0+0.3,1.5+0.15);
\draw[cyan] node at (1.55,1.8) {\scriptsize{$6,c$}};
\draw [line width=0.9mm,dashed] (0.7,1.5) -- (1.3,1.5);
\end{scope}

\begin{scope}[shift={(-1.33,-0.65)}]
\draw[cyan] (1.0,1.5) -- (1.0+0.3,1.5+0.15);
\draw[cyan] (1.0+0.3,1.5+0.15) -- (1.0+0.3,1.5-0.15);
\draw[cyan] (1.0+0.3,1.5-0.15) -- (1.0,1.5);
\draw[cyan] node at (1.6,1.75) {\scriptsize{$7,d$}};
\draw [line width=0.9mm,dashed] (1.35,1.5) -- (1.35+0.6,1.5);
\end{scope}

\end{scope}

\begin{scope}[shift={(+2.8*2,0.0)}]
\draw[black] node at (1.0,2.15) {\scriptsize{$\mathcal S_{\text{eq}}^{(3)}$}};
\draw[black,thick] (0,0) rectangle (1.8,1.8);
\draw[red] (0.3,0.5) -- (0.3+0.3,0.5);
\draw[red] (0.3+0.3,0.5) -- (0.3+0.15,0.5/2);
\draw[red] (0.3+0.15,0.5/2) -- (0.3,0.5);

\draw[red] (0.3,0.5) -- (0.3+0.3,0.5);
\draw[red] (0.3+0.3,0.5) -- (0.3+0.15,0.5+0.25);
\draw[red] (0.3+0.15,0.5+0.25) -- (0.3,0.5);
\draw[red] node at (0.3+0.6,0.5) {\scriptsize{$8,a$}};

\draw[blue] (1.0,1.5) -- (1.0+0.3,1.5+0.15);
\draw[blue] (1.0+0.3,1.5+0.15) -- (1.0+0.3,1.5-0.15);
\draw[blue] (1.0+0.3,1.5-0.15) -- (1.0,1.5);

\draw[blue] (1.0+0.3,1.5+0.15) -- (1.0+0.3,1.5-0.15);
\draw[blue] (1.0+0.3,1.5-0.15) -- (1.6,1.5);
\draw[blue] (1.6,1.5) -- (1.0+0.3,1.5+0.15);
\draw[blue] node at (1.0+0.3,1.5-0.25) {\scriptsize{$9,b$}};

\begin{scope}[shift={(-1.33,-0.65)}]
\draw[cyan] (1.0,1.5) -- (1.0+0.3,1.5+0.15);
\draw[cyan] (1.0+0.3,1.5+0.15) -- (1.0+0.3,1.5-0.15);
\draw[cyan] (1.0+0.3,1.5-0.15) -- (1.0,1.5);
\draw[cyan] node at (1.65,1.75) {\scriptsize{$11,d$}};
\draw [line width=0.9mm,dashed] (1.35,1.5) -- (1.35+0.6,1.5);
\end{scope}

\end{scope}

\end{tikzpicture}}}
\hfill \null
\\
\null \hfill
\subfloat[\label{fig:conductor_traverse_meshing3} ] {\resizebox{1.03\columnwidth} {!} {\begin{tikzpicture}[scale=1]

\draw[black] node at (1.0,2.15) {\scriptsize{$\mathcal S_{\text{eq}}^{(1)}$}};
\draw[black,thick] (0,0) rectangle (1.8,1.8);
\draw[red] (0.3,0.5) -- (0.3+0.3,0.5);
\draw[red] (0.3+0.3,0.5) -- (0.3+0.15,0.5/2);
\draw[red] (0.3+0.15,0.5/2) -- (0.3,0.5);

\draw[red] (0.3,0.5) -- (0.3+0.3,0.5);
\draw[red] (0.3+0.3,0.5) -- (0.3+0.15,0.5+0.25);
\draw[red] (0.3+0.15,0.5+0.25) -- (0.3,0.5);
\draw[red] node at (0.3+0.6,0.5) {\scriptsize{$1,a$}};

\draw[blue] (1.0,1.5) -- (1.0+0.3,1.5+0.15);
\draw[blue] (1.0+0.3,1.5+0.15) -- (1.0+0.3,1.5-0.15);
\draw[blue] (1.0+0.3,1.5-0.15) -- (1.0,1.5);

\draw[blue] (1.0+0.3,1.5+0.15) -- (1.0+0.3,1.5-0.15);
\draw[blue] (1.0+0.3,1.5-0.15) -- (1.6,1.5);
\draw[blue] (1.6,1.5) -- (1.0+0.3,1.5+0.15);
\draw[blue] node at (1.0+0.3,1.5-0.28) {\scriptsize{$2,b$}};

\begin{scope}[shift={(0.53,-0.65)}]
\draw[cyan] (1.0+0.3,1.5+0.15) -- (1.0+0.3,1.5-0.15);
\draw[cyan] (1.0+0.3,1.5-0.15) -- (1.6,1.5);
\draw[cyan] (1.6,1.5) -- (1.0+0.3,1.5+0.15);
\draw[cyan] node at (1.55,1.8) {\scriptsize{$3,c$}};
\draw [line width=0.9mm,dashed] (0.7,1.5) -- (1.3,1.5);
\end{scope}

\begin{scope}[shift={(-1.33,-0.65+0.45)}]
\draw[cyan] (1.0,1.5) -- (1.0+0.3,1.5+0.15);
\draw[cyan] (1.0+0.3,1.5+0.15) -- (1.0+0.3,1.5-0.15);
\draw[cyan] (1.0+0.3,1.5-0.15) -- (1.0,1.5);
\draw[cyan] node at (1.6,1.75) {\scriptsize{$4,d$}};
\draw [line width=0.9mm,dashed] (1.35,1.5) -- (1.35+0.6,1.5);
\end{scope}

\begin{scope}[shift={(+2.8,0.0)}]
\draw[black] node at (1.0,2.15) {\scriptsize{$\mathcal S_{\text{eq}}^{(2)}$}};
\draw[black,thick] (0,0) rectangle (1.8,1.8);
\draw[red] (0.3,0.5) -- (0.3+0.3,0.5);
\draw[red] (0.3+0.3,0.5) -- (0.3+0.15,0.5/2);
\draw[red] (0.3+0.15,0.5/2) -- (0.3,0.5);

\draw[red] (0.3,0.5) -- (0.3+0.3,0.5);
\draw[red] (0.3+0.3,0.5) -- (0.3+0.15,0.5+0.25);
\draw[red] (0.3+0.15,0.5+0.25) -- (0.3,0.5);
\draw[red] node at (0.3+0.6,0.5) {\scriptsize{$5,a$}};

\draw[blue] (1.0,1.5) -- (1.0+0.3,1.5+0.15);
\draw[blue] (1.0+0.3,1.5+0.15) -- (1.0+0.3,1.5-0.15);
\draw[blue] (1.0+0.3,1.5-0.15) -- (1.0,1.5);

\draw[blue] (1.0+0.3,1.5+0.15) -- (1.0+0.3,1.5-0.15);
\draw[blue] (1.0+0.3,1.5-0.15) -- (1.6,1.5);
\draw[blue] (1.6,1.5) -- (1.0+0.3,1.5+0.15);
\draw[blue] node at (1.0+0.3-.4,1.5-0.28+.4) {\scriptsize{$6,b$}};

\begin{scope}[shift={(0.53,-0.65+0.2)}]
\draw[cyan] (1.0+0.3,1.5+0.15) -- (1.0+0.3,1.5-0.15);
\draw[cyan] (1.0+0.3,1.5-0.15) -- (1.6,1.5);
\draw[cyan] (1.6,1.5) -- (1.0+0.3,1.5+0.15);
\draw[cyan] node at (1.55,1.8) {\scriptsize{$7,c$}};
\draw [line width=0.9mm,dashed] (0.7,1.5) -- (1.3,1.5);
\end{scope}

\begin{scope}[shift={(-1.33,-0.65)}]
\draw[cyan] (1.0,1.5) -- (1.0+0.3,1.5+0.15);
\draw[cyan] (1.0+0.3,1.5+0.15) -- (1.0+0.3,1.5-0.15);
\draw[cyan] (1.0+0.3,1.5-0.15) -- (1.0,1.5);
\draw[cyan] node at (1.6,1.75) {\scriptsize{$8,d$}};
\draw [line width=0.9mm,dashed] (1.35,1.5) -- (1.35+0.6,1.5);
\end{scope}

\end{scope}

\begin{scope}[shift={(+2.8*2,0.0)}]
\draw[black] node at (1.0,2.15) {\scriptsize{$\mathcal S_{\text{eq}}^{(3)}$}};
\draw[black,thick] (0,0) rectangle (1.8,1.8);
\draw[red] (0.3,0.5) -- (0.3+0.3,0.5);
\draw[red] (0.3+0.3,0.5) -- (0.3+0.15,0.5/2);
\draw[red] (0.3+0.15,0.5/2) -- (0.3,0.5);

\draw[red] (0.3,0.5) -- (0.3+0.3,0.5);
\draw[red] (0.3+0.3,0.5) -- (0.3+0.15,0.5+0.25);
\draw[red] (0.3+0.15,0.5+0.25) -- (0.3,0.5);
\draw[red] node at (0.3+0.6,0.5) {\scriptsize{$9,a$}};

\draw[blue] (1.0,1.5) -- (1.0+0.3,1.5+0.15);
\draw[blue] (1.0+0.3,1.5+0.15) -- (1.0+0.3,1.5-0.15);
\draw[blue] (1.0+0.3,1.5-0.15) -- (1.0,1.5);

\draw[blue] (1.0+0.3,1.5+0.15) -- (1.0+0.3,1.5-0.15);
\draw[blue] (1.0+0.3,1.5-0.15) -- (1.6,1.5);
\draw[blue] (1.6,1.5) -- (1.0+0.3,1.5+0.15);
\draw[blue] node at (1.0+0.3,1.5-0.28) {\scriptsize{$10,b$}};

\begin{scope}[shift={(0.53,-0.65)}]
\draw[cyan] (1.0+0.3,1.5+0.15) -- (1.0+0.3,1.5-0.15);
\draw[cyan] (1.0+0.3,1.5-0.15) -- (1.6,1.5);
\draw[cyan] (1.6,1.5) -- (1.0+0.3,1.5+0.15);
\draw[cyan] node at (1.6,1.8) {\scriptsize{$11,c$}};
\draw [line width=0.9mm,dashed] (0.7,1.5) -- (1.3,1.5);
\end{scope}

\begin{scope}[shift={(-1.33,-0.65+0.2)}]
\draw[cyan] (1.0,1.5) -- (1.0+0.3,1.5+0.15);
\draw[cyan] (1.0+0.3,1.5+0.15) -- (1.0+0.3,1.5-0.15);
\draw[cyan] (1.0+0.3,1.5-0.15) -- (1.0,1.5);
\draw[cyan] node at (1.65,1.75) {\scriptsize{$12,d$}};
\draw [line width=0.9mm,dashed] (1.35,1.5) -- (1.35+0.6,1.5);
\end{scope}

\end{scope}

\end{tikzpicture}}}
\hfill \null
\caption{Top view of three sample arrays of three equivalent surfaces with their  basis functions. In reality, the space between equivalent surfaces is zero. The full RWG basis functions are shown in blue and red and the half RWG basis functions are shown in cyan. Each basis function has a local and a global identification number. Local identification numbers are denoted by $a, b, c, d$ and global identification numbers are denoted by $1, 2, \dots, 12$. (a) First scenario: a structure where PEC traces traverse the equivalent surfaces in the same location in all unit cells. (b) Second scenario: a structure where PEC traces traverse equivalent surfaces in the same location, but conductor traverse is only present between some unit cells. (c) Third scenario: a structure where PEC traces traverse the equivalent surfaces in different and arbitrary locations. Dashed lines show the position where the PEC traces traverse the equivalent surfaces.}
\label{fig:array_identical_meshing}
\end{figure}
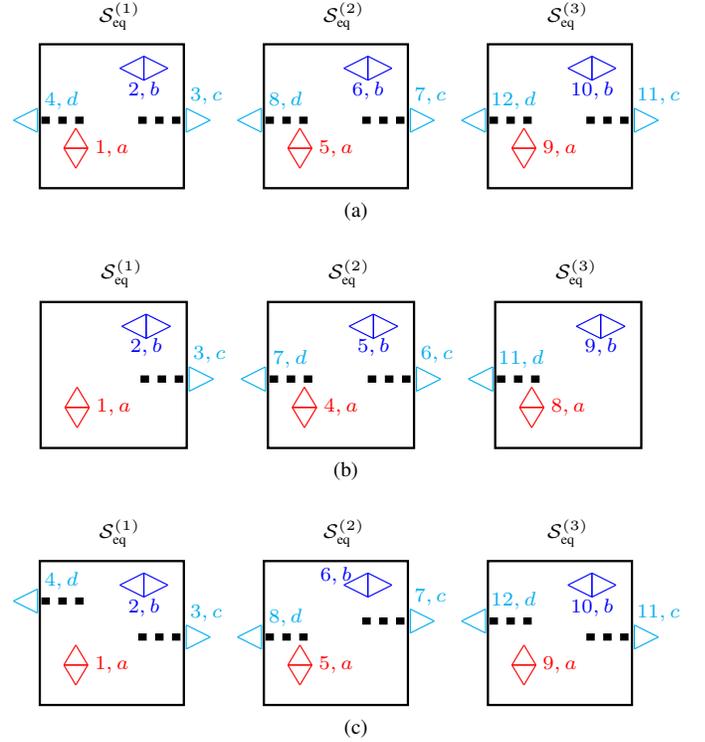

In order to cast $\textbf Z_0$ in~\eqref{U0tZU0X_matrix-1} into a Toeplitz matrix~\cite{FFT1}, which is necessary for the FFT-based acceleration technique, we must have the following conditions:

\begin{enumerate}[label=(\roman*)]
\item The equivalent surfaces must be all identical;
\item 
The meshes on the equivalent surfaces must be all identical; 
\item 
The basis functions on the equivalent surfaces must all be identical. 
\end{enumerate}
In a structure with no conductors traversing the equivalent surfaces, these three conditions can be easily met by using the same mesh on each equivalent surface $\mathcal S_{\text{eq}}^{(m)}$~\cite{Utkarsh_macromodel_2020}. Similar to the previous work, the equivalent surfaces in the proposed approach are all identical while the unit cells in the EM surfaces can be non-identical. Hence, the condition (i) is met here. In presence of conductor traverse, three scenarios can arise. In order to carefully address these scenarios, let us consider three structures each with different positioning of conductors traversing the equivalent surfaces. Hence, three sample arrays of three uniformly spaced equivalent surfaces are depicted in Fig.~\ref{fig:array_identical_meshing}. If conductors traverse the equivalent surfaces in the same location in all unit cells, as in Fig.~\ref{fig:array_identical_meshing}\subref{fig:conductor_traverse_meshing1}, the same meshes and basis functions can be used on all side faces of the equivalent surfaces. The three conditions that enable FFT acceleration are therefore met. This first scenario arises in some practical EM surfaces, for example when the conductor traversing the equivalent surfaces is a slotted ground plane, and slots are fully contained inside each unit cell~\cite{EMS_Connected5, EMS_Connected8}. The second scenario arises when conductors traverse equivalent surfaces in the same location, but conductor traverse is only present between some unit cells~\cite{EMS_Connected9}, as shown in Fig.~\ref{fig:array_identical_meshing}\subref{fig:conductor_traverse_meshing2}. In this case, to maintain the periodicity of the basis functions (condition (iii)), the half RWG basis functions need to be introduced on the other unit cells even where there is no conductor traverse, i.e. edges on the left side of $\mathcal S_{\text{eq}}^{(1)}$ and the right side of $\mathcal S_{\text{eq}}^{(3)}$ need to be expanded by the half RWG basis functions (in Fig.~\ref{fig:array_identical_meshing}\subref{fig:conductor_traverse_meshing2}). Hence, this na\"ive modification enable the FFT acceleration for simulation of this structure. The third scenario is the case where conductor traverse happens in different and arbitrary locations~\cite{EMS_Connected7}, as shown in Fig.~\ref{fig:array_identical_meshing}\subref{fig:conductor_traverse_meshing3}. To apply FFT acceleration to this latter case, one needs to generate a mesh and set of basis functions for the side surfaces of all unit cells that satisfy conditions (ii) and (iii) defined above. Although this is conceivable, a simpler approach may be to resort to other acceleration techniques which do not require a periodic mesh, such as the AIM or the MLFMM. Comparing FFT-based acceleration with these more general alternatives will be the subject of future work. It is worth mentioning that FFT-based acceleration technique is simpler to implement than AIM and MLFMM.

\section{Numerical Results}
\label{Numerical Results}   
In order to validate the proposed method for simulation of EM surfaces, three examples are presented. The results are compared against the FEKO~\cite{FEKO} and ANSYS HFSS~\cite{HFSS} commercial solvers. The numerical results are obtained on a single Intel Xeon Silver $4110$ processor running with a single core at $2.10$ GHz and $512$ GB of RAM.

\subsection{Array of Rectangular PEC Patches}
\label{Example1_Array_Rectangular_Patches}
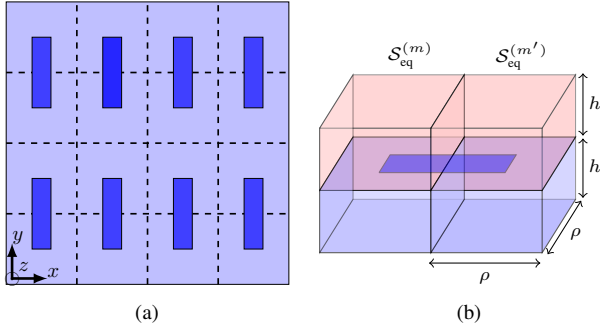
\begin{figure}
\hfill \null
\\
\null \hfill
\subfloat[\label{fig:array_rectangular_patches_a} ] 
{\resizebox{0.45\columnwidth} {!} {\begin{tikzpicture}[scale=1]

\filldraw[fill=blue!25!white, draw=black] (0.0,0.0) rectangle (1.2*4,1.2*4);
\draw[black,thick,dashed] (1.2,0.0) -- (1.2,1.2*4);
\draw[black,thick,dashed] (1.2*2,0.0) -- (1.2*2,1.2*4);
\draw[black,thick,dashed] (1.2*3,0.0) -- (1.2*3,1.2*4);

\draw[black,thick,dashed] (0.0,1.2) -- (1.2*4,1.2);
\draw[black,thick,dashed] (0.0,1.2*2) -- (1.2*4,1.2*2);
\draw[black,thick,dashed] (0.0,1.2*3) -- (1.2*4,1.2*3);

\begin{scope}[shift={(1.3,0.1)}]
\draw node at (-1.1,.7) {$y$};
\draw [black][line width=0.45mm,-latex] (-1.2,0.0) -- (-1.2,.6);
\draw node at (-0.5,.1) {$x$};
\draw node at (-1.2,0.0) {$\odot$};
\draw node at (-1.0,.2) {$z$};
\draw [black][line width=0.45mm,-latex] (-1.2,0.0) -- (-.6,0.0);
\end{scope}

\filldraw[fill=blue!75!white, draw=black] (0.6-.16,0.6) rectangle (0.6+.16,.6+1.2);
\begin{scope}[shift={(1.2,0.0)}]
\filldraw[fill=blue!75!white, draw=black] (0.6-.16,0.6) rectangle (0.6+.16,.6+1.2);
\end{scope}
\begin{scope}[shift={(1.2*2,0.0)}]
\filldraw[fill=blue!75!white, draw=black] (0.6-.16,0.6) rectangle (0.6+.16,.6+1.2);
\end{scope}
\begin{scope}[shift={(1.2*3,0.0)}]
\filldraw[fill=blue!75!white, draw=black] (0.6-.16,0.6) rectangle (0.6+.16,.6+1.2);
\end{scope}

\begin{scope}[shift={(0.0,1.2*2)}]
\filldraw[fill=blue!75!white, draw=black] (0.6-.16,0.6) rectangle (0.6+.16,.6+1.2);
\begin{scope}[shift={(1.2,0.0)}]
\filldraw[fill=blue!85!white, draw=black] (0.6-.16,0.6) rectangle (0.6+.16,.6+1.2);
\end{scope}
\begin{scope}[shift={(1.2*2,0.0)}]
\filldraw[fill=blue!75!white, draw=black] (0.6-.16,0.6) rectangle (0.6+.16,.6+1.2);
\end{scope}
\begin{scope}[shift={(1.2*3,0.0)}]
\filldraw[fill=blue!75!white, draw=black] (0.6-.16,0.6) rectangle (0.6+.16,.6+1.2);
\end{scope}
\end{scope}

\end{tikzpicture}}}
\hfill 
\subfloat[\label{fig:sample_cell_equivalent_b}
]{\resizebox{0.45\columnwidth} {!}  {\begin{tikzpicture}[xscale= 0.25, yscale =0.4]
\pgfmathsetmacro{\h}{0.7}

\begin{scope}[shift = {($(0,-1.5*\h)$)}, every node/.append style={xslant=1},xslant=1, xscale= 1, yscale = 0.3]
\coordinate (p1) at (-3.0,-3.0);
\coordinate (p2) at (3.0,-3.0);
\coordinate (p3) at (3.0, 3.0);
\coordinate (p4) at (-3.0,3.0);
\draw[fill=red!40!white,opacity=0.4] (p1)--(p2)--(p3)--(p4)--(p1);
\end{scope}

\draw[fill=red!40!white, opacity = 0.4] (p1) --++ (0,3*\h) --++ (6.0,0) -- (p2) -- (p1);

\begin{scope}[shift = {($(0,1.5*\h)$)}, every node/.append style={xslant=1},xslant=1, xscale= 1, yscale = 0.3]
\coordinate (p1) at (-3.0,-3.0);
\coordinate (p2) at (3.0,-3.0);
\coordinate (p3) at (3.0, 3.0);
\coordinate (p4) at (-3.0,3.0);
\draw[fill=red!40!white,opacity=0.4] (p1)--(p2)--(p3)--(p4)--(p1);
\end{scope}

\draw[fill=red!40!white, opacity = 0.4] (p2) -- (p3) -- ($(p3) + (0,-3*\h)$) -- ($(p2) + (0,-3*\h)$) -- (p2);

\draw[fill=red!10!white, opacity = 0.4] (p1) -- (p4) -- ($(p4) + (0,-3*\h)$) -- ($(p1) + (0,-3*\h)$) -- (p1);

\begin{scope}[shift={(0.0,-3*\h)}]
\begin{scope}[shift = {($(0,-1.5*\h)$)}, every node/.append style={xslant=1},xslant=1, xscale= 1, yscale = 0.3]
\coordinate (p1) at (-3.0,-3.0);
\coordinate (p2) at (3.0,-3.0);
\coordinate (p3) at (3.0, 3.0);
\coordinate (p4) at (-3.0,3.0);
\draw[fill=blue!35!white,opacity=0.4] (p1)--(p2)--(p3)--(p4)--(p1);
\end{scope}

\draw[fill=blue!35!white, opacity = 0.4] (p1) --++ (0,3*\h) --++ (6.0,0) -- (p2) -- (p1);

\begin{scope}[shift = {($(0,1.5*\h)$)}, every node/.append style={xslant=1},xslant=1, xscale= 1, yscale = 0.3]
\coordinate (p1) at (-3.0,-3.0);
\coordinate (p2) at (3.0,-3.0);
\coordinate (p3) at (3.0, 3.0);
\coordinate (p4) at (-3.0,3.0);
\draw[fill=blue!35!white,opacity=0.4] (p1)--(p2)--(p3)--(p4)--(p1);
\end{scope}

\draw[fill=blue!35!white, opacity = 0.4] (p2) -- (p3) -- ($(p3) + (0,-3*\h)$) -- ($(p2) + (0,-3*\h)$) -- (p2);

\draw[fill=blue!35!white, opacity = 0.4] (p1) -- (p4) -- ($(p4) + (0,-3*\h)$) -- ($(p1) + (0,-3*\h)$) -- (p1);

\end{scope}

\begin{scope}[shift={(6.0,0.0)}]
\begin{scope}[shift = {($(0,-1.5*\h)$)}, every node/.append style={xslant=1},xslant=1, xscale= 1, yscale = 0.3]
\coordinate (p1) at (-3.0,-3.0);
\coordinate (p2) at (3.0,-3.0);
\coordinate (p3) at (3.0, 3.0);
\coordinate (p4) at (-3.0,3.0);
\draw[fill=red!40!white,opacity=0.4] (p1)--(p2)--(p3)--(p4)--(p1);
\end{scope}

\draw[fill=red!40!white, opacity = 0.4] (p1) --++ (0,3*\h) --++ (6.0,0) -- (p2) -- (p1);

\begin{scope}[shift = {($(0,1.5*\h)$)}, every node/.append style={xslant=1},xslant=1, xscale= 1, yscale = 0.3]
\coordinate (p1) at (-3.0,-3.0);
\coordinate (p2) at (3.0,-3.0);
\coordinate (p3) at (3.0, 3.0);
\coordinate (p4) at (-3.0,3.0);
\draw[fill=red!40!white,opacity=0.4] (p1)--(p2)--(p3)--(p4)--(p1);
\end{scope}

\draw[fill=red!40!white, opacity = 0.4] (p2) -- (p3) -- ($(p3) + (0,-3*\h)$) -- ($(p2) + (0,-3*\h)$) -- (p2);

\draw[fill=red!10!white, opacity = 0.4] (p1) -- (p4) -- ($(p4) + (0,-3*\h)$) -- ($(p1) + (0,-3*\h)$) -- (p1);

\begin{scope}[shift={(0.0,-3*\h)}]
\begin{scope}[shift = {($(0,-1.5*\h)$)}, every node/.append style={xslant=1},xslant=1, xscale= 1, yscale = 0.3]
\coordinate (p1) at (-3.0,-3.0);
\coordinate (p2) at (3.0,-3.0);
\coordinate (p3) at (3.0, 3.0);
\coordinate (p4) at (-3.0,3.0);
\draw[fill=blue!40!white,opacity=0.4] (p1)--(p2)--(p3)--(p4)--(p1);
\end{scope}

\draw[fill=blue!40!white, opacity = 0.4] (p1) --++ (0,3*\h) --++ (6.0,0) -- (p2) -- (p1);

\begin{scope}[shift = {($(0,1.5*\h)$)}, every node/.append style={xslant=1},xslant=1, xscale= 1, yscale = 0.3]
\coordinate (p1) at (-3.0,-3.0);
\coordinate (p2) at (3.0,-3.0);
\coordinate (p3) at (3.0, 3.0);
\coordinate (p4) at (-3.0,3.0);
\draw[fill=blue!40!white,opacity=0.4] (p1)--(p2)--(p3)--(p4)--(p1);
\end{scope}

\draw[fill=blue!40!white, opacity = 0.4] (p2) -- (p3) -- ($(p3) + (0,-3*\h)$) -- ($(p2) + (0,-3*\h)$) -- (p2);

\draw[fill=blue!20!white, opacity = 0.4] (p1) -- (p4) -- ($(p4) + (0,-3*\h)$) -- ($(p1) + (0,-3*\h)$) -- (p1);

\begin{scope}[shift = {($(-6.0,1.5*\h)$)}, every node/.append style={xslant=1},xslant=1, xscale= 1, yscale = 0.3]
\coordinate (p1) at (-0.4,-1.0);
\coordinate (p2) at (2*3.2,-1.0);
\coordinate (p3) at (2*3.2, 1.0);
\coordinate (p4) at (-0.4,1.0);
\draw[fill=blue!99!white,opacity=0.4] (p1)--(p2)--(p3)--(p4)--(p1);
\end{scope}

\end{scope}

\end{scope}

\draw node at (1.0,2.6) {\scriptsize $\mathcal S_{\text{eq}}^{(m)}$};
\draw node at (7.0,2.6) {\scriptsize $\mathcal S_{\text{eq}}^{(m')}$};

\draw[black,<->] (10.3,-0.1) -- (10.3,2.0);
\draw node at (10.95,1.0) {\scriptsize $h$};
\draw[black,<->] (10.3,-2.2) -- (10.3,-.2);\\
\draw node at (10.95,-1.1) {\scriptsize $h$};

\draw[black,<->] (8.4,-4.1) -- (10.2,-2.3);\\
\draw node at (9.9,-3.4) {\scriptsize $\rho$};
\draw[black,<->] (2.0,-4.4) -- (8.0,-4.4);
\draw node at (5.0,-4.9) {\scriptsize $\rho$};

\end{tikzpicture}
}}
\hfill \null
\caption{(a) Top view of $4 \times 2$ array of rectangular PEC patches. The dielectric layer on the top layer is eliminated for better visualization of array distribution. Dashed lines represent the boundary of the equivalent surfaces in the proposed approach which generate an array of $4 \times 4$ unit cells. It is clear that the PEC patches traverse the macromodel equivalent surfaces. (b) Two adjacent unit cells of the array of rectangular patches in the proposed solver with $\rho =13.5$ mm and $h=1$ mm.} 
\label{fig:Array_Rectangular_Patches}
\end{figure}

\begin{figure}
\hfill \null
\\
\null \hfill
\subfloat[\label{fig:array2_rectangular_patches_a} ] 
{\resizebox{0.45\columnwidth} {!} {\begin{tikzpicture}[scale=1]

\filldraw[fill=blue!25!white, draw=black] (0.0,0.0) rectangle (1.2*4,1.2*4);
\draw[black,thick,dashed] (1.2,0.0) -- (1.2,1.2*4);
\draw[black,thick,dashed] (1.2*2,0.0) -- (1.2*2,1.2*4);
\draw[black,thick,dashed] (1.2*3,0.0) -- (1.2*3,1.2*4);

\draw[black,thick,dashed] (0.0,1.2*2) -- (1.2*4,1.2*2);

\begin{scope}[shift={(1.3,0.1)}]
\draw node at (-1.1,.7) {$y$};
\draw [black][line width=0.45mm,-latex] (-1.2,0.0) -- (-1.2,.6);
\draw node at (-0.5,.1) {$x$};
\draw node at (-1.2,0.0) {$\odot$};
\draw node at (-1.0,.2) {$z$};
\draw [black][line width=0.45mm,-latex] (-1.2,0.0) -- (-.6,0.0);
\end{scope}

\filldraw[fill=blue!75!white, draw=black] (0.6-.16,0.6) rectangle (0.6+.16,.6+1.2);
\begin{scope}[shift={(1.2,0.0)}]
\filldraw[fill=blue!75!white, draw=black] (0.6-.16,0.6) rectangle (0.6+.16,.6+1.2);
\end{scope}
\begin{scope}[shift={(1.2*2,0.0)}]
\filldraw[fill=blue!75!white, draw=black] (0.6-.16,0.6) rectangle (0.6+.16,.6+1.2);
\end{scope}
\begin{scope}[shift={(1.2*3,0.0)}]
\filldraw[fill=blue!75!white, draw=black] (0.6-.16,0.6) rectangle (0.6+.16,.6+1.2);
\end{scope}

\begin{scope}[shift={(0.0,1.2*2)}]
\filldraw[fill=blue!75!white, draw=black] (0.6-.16,0.6) rectangle (0.6+.16,.6+1.2);
\begin{scope}[shift={(1.2,0.0)}]
\filldraw[fill=blue!85!white, draw=black] (0.6-.16,0.6) rectangle (0.6+.16,.6+1.2);
\end{scope}
\begin{scope}[shift={(1.2*2,0.0)}]
\filldraw[fill=blue!75!white, draw=black] (0.6-.16,0.6) rectangle (0.6+.16,.6+1.2);
\end{scope}
\begin{scope}[shift={(1.2*3,0.0)}]
\filldraw[fill=blue!75!white, draw=black] (0.6-.16,0.6) rectangle (0.6+.16,.6+1.2);
\end{scope}
\end{scope}

\end{tikzpicture}}}
\hfill 
\subfloat[\label{fig:sample_cell2_equivalent_b}
]{\resizebox{0.45\columnwidth} {!}  {\begin{tikzpicture}[xscale= 0.25, yscale =0.4]
\pgfmathsetmacro{\h}{0.7}

\begin{scope}[shift = {($(0,-1.5*\h)$)}, every node/.append style={xslant=1},xslant=1, xscale= 1, yscale = 0.3]
\coordinate (p1) at (-3.0*2,-3.0);
\coordinate (p2) at (3.0*2,-3.0);
\coordinate (p3) at (3.0*2, 3.0);
\coordinate (p4) at (-3.0*2,3.0);
\draw[fill=red!40!white,opacity=0.4] (p1)--(p2)--(p3)--(p4)--(p1);
\end{scope}

\draw[fill=red!40!white, opacity = 0.4] (p1) --++ (0,3*\h) --++ (6.0*2,0) -- (p2) -- (p1);

\begin{scope}[shift = {($(0,1.5*\h)$)}, every node/.append style={xslant=1},xslant=1, xscale= 1, yscale = 0.3]
\coordinate (p1) at (-3.0*2,-3.0);
\coordinate (p2) at (3.0*2,-3.0);
\coordinate (p3) at (3.0*2,3.0);
\coordinate (p4) at (-3.0*2,3.0);
\draw[fill=red!40!white,opacity=0.4] (p1)--(p2)--(p3)--(p4)--(p1);
\end{scope}

\draw[fill=red!40!white, opacity = 0.4] (p2) -- (p3) -- ($(p3) + (0,-3*\h)$) -- ($(p2) + (0,-3*\h)$) -- (p2);

\draw[fill=red!10!white, opacity = 0.4] (p1) -- (p4) -- ($(p4) + (0,-3*\h)$) -- ($(p1) + (0,-3*\h)$) -- (p1);

\begin{scope}[shift={(0.0,-3*\h)}]
\begin{scope}[shift = {($(0,-1.5*\h)$)}, every node/.append style={xslant=1},xslant=1, xscale= 1, yscale = 0.3]
\coordinate (p1) at (-3.0*2,-3.0);
\coordinate (p2) at (3.0*2,-3.0);
\coordinate (p3) at (3.0*2, 3.0);
\coordinate (p4) at (-3.0*2,3.0);
\draw[fill=blue!35!white,opacity=0.4] (p1)--(p2)--(p3)--(p4)--(p1);
\end{scope}

\draw[fill=blue!35!white, opacity = 0.4] (p1) --++ (0,3*\h) --++ (6.0*2,0) -- (p2) -- (p1);

\begin{scope}[shift = {($(0,1.5*\h)$)}, every node/.append style={xslant=1},xslant=1, xscale= 1, yscale = 0.3]
\coordinate (p1) at (-3.0*2,-3.0);
\coordinate (p2) at (3.0*2,-3.0);
\coordinate (p3) at (3.0*2, 3.0);
\coordinate (p4) at (-3.0*2,3.0);
\draw[fill=blue!35!white,opacity=0.4] (p1)--(p2)--(p3)--(p4)--(p1);
\end{scope}

\draw[fill=blue!35!white, opacity = 0.4] (p2) -- (p3) -- ($(p3) + (0,-3*\h)$) -- ($(p2) + (0,-3*\h)$) -- (p2);

\draw[fill=blue!35!white, opacity = 0.4] (p1) -- (p4) -- ($(p4) + (0,-3*\h)$) -- ($(p1) + (0,-3*\h)$) -- (p1);

\end{scope}

\begin{scope}[shift = {($(-3.0,-1.5*\h)$)}, every node/.append style={xslant=1},xslant=1, xscale= 1, yscale = 0.3]
\coordinate (p1) at (-0.4,-1.0);
\coordinate (p2) at (2*3.2,-1.0);
\coordinate (p3) at (2*3.2, 1.0);
\coordinate (p4) at (-0.4,1.0);
\draw[fill=blue!99!white,opacity=0.4] (p1)--(p2)--(p3)--(p4)--(p1);
\end{scope}

\draw node at (1.0,2.6) {\scriptsize $\mathcal S_{\text{eq}}^{(m)}$};

\begin{scope}[shift={(-3.0,0.0)}]
\draw[black,<->] (10.3,-0.1) -- (10.3,2.0);
\draw node at (10.95,1.0) {\scriptsize $h$};
\draw[black,<->] (10.3,-2.2) -- (10.3,-.2);\\
\draw node at (10.95,-1.1) {\scriptsize $h$};

\draw[black,<->] (8.4,-4.1) -- (10.2,-2.3);\\
\draw node at (9.9,-3.4) {\scriptsize $\rho$};
\draw[black,<->] (2.0-6.0,-4.4) -- (8.0,-4.4);
\draw node at (2.5,-4.9) {\scriptsize $2\rho$};
\end{scope}

\end{tikzpicture}
}}
\hfill \null
\caption{(a) The same structure of array of rectangular PEC patches in Fig.~\ref{fig:Array_Rectangular_Patches} considered for macromodeling approach in the previous work~\cite{Utkarsh_macromodel_2020}. Dashed lines represent the larger equivalent surfaces used in the previous work to avoid cutting the PEC patches. This macromodeling strategy produces an array of $4 \times 2$ unit cells. (b) A unit cell of the array of rectangular patches which is not allowed to traverse the equivalent surfaces. Here, $\rho =13.5$ mm and $h=1$ mm.}
\label{fig:Array_Rectangular_Patches2}
\end{figure}
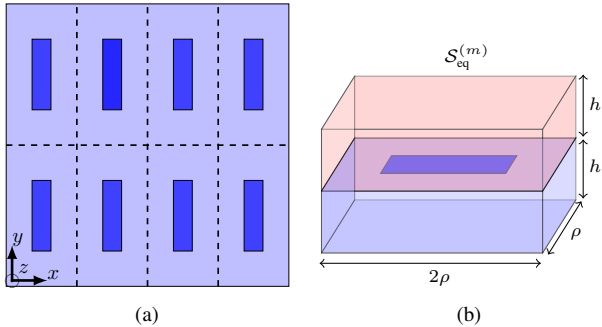

For the first example, a simple test is considered to validate the accuracy of the proposed approach. Hence, we assume a scattering problem involving an array of $4 \times 2$ rectangular PEC patch elements, as shown in Fig.~\ref{fig:Array_Rectangular_Patches}. The dashed lines in Fig.~\ref{fig:Array_Rectangular_Patches}\subref{fig:array_rectangular_patches_a} represent the boundaries of the equivalent surfaces used in the proposed approach, which generate an array of $4 \times 4$ unit cells. As shown in Fig.~\ref{fig:Array_Rectangular_Patches}\subref{fig:sample_cell_equivalent_b}, the PEC patches traverse the macromodeling boundaries. We choose this particular example since it can also be solved by the macromodeling approach developed in the previous work~\cite{Utkarsh_macromodel_2020} by changing the equivalence surfaces.
As reviewed earlier, the PEC patches were not allowed to traverse the equivalent surfaces in the previous work. Hence, as depicted in Fig.~\ref{fig:Array_Rectangular_Patches2}\subref{fig:array2_rectangular_patches_a}, the equivalent surfaces represented by dashed lines need to be larger to avoid cutting the PEC patches. This strategy of macromodeling produces an array of $4 \times 2$ unit cells and is used in this section to validate the accuracy of the current modeling of the proposed approach.  

In this example, the array substrate has two layers, each with a thickness of $1.0$ mm. As shown in Fig.~\ref{fig:Array_Rectangular_Patches}\subref{fig:sample_cell_equivalent_b}, each unit cell has dimensions of $13.5$ mm $\times~ 13.5$ mm $\times~2.0$ mm with PEC patch size of $6.75$ mm $\times~2.7$ mm. The equivalent surface $\mathcal S_{\text{eq}}^{(m)}$ encapsulating each unit cell has dimensions of $13.5$ mm $\times~13.5~\times$ $2.0$ mm. The bottom substrate layer has relative permittivity of $\varepsilon_r = 2.2$ and the top substrate layer has relative permittivity of $\varepsilon_r = 3.0$. The incident field is produced by an electric dipole with moment $I\cdot l = 1$ A$\cdot$m , which operates at $9.6$ GHz. The dipole is located at $(0,0,20\lambda_0=0.625~\text{m})$ and is directed along the $x$-axis. The rectangular PEC patch is discretized with triangular meshes with a characteristic length of $1.2$ mm. The characteristic mesh length along the dielectric substrate is chosen to be $2.0$ mm. The same discretization sample is chosen for the FEKO-MLFMM simulation. In the proposed approach, the equivalent surface of each unit cell is discretized with a characteristic length of $2.5$ mm, which is chosen based on experiments. 

\begin{figure}[t]
\begin{center}
\subfloat[\label{fig:Directivity_ex1_phi_0} ] 
{\resizebox{0.9\columnwidth} {!}
{\includegraphics{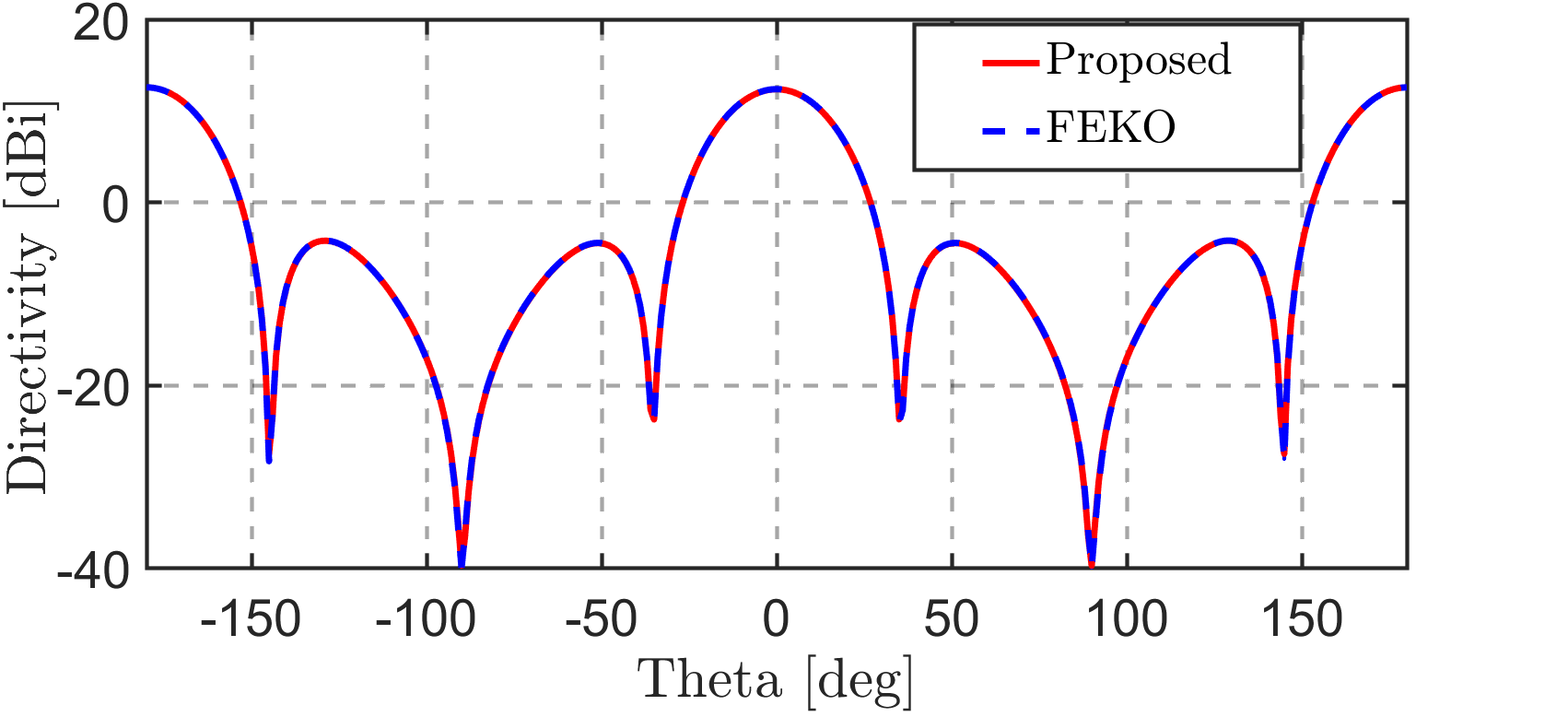}}} 
\\
\subfloat[\label{fig:Directivity_ex1_phi_0} ] 
{\resizebox{0.9\columnwidth} {!}
{\includegraphics{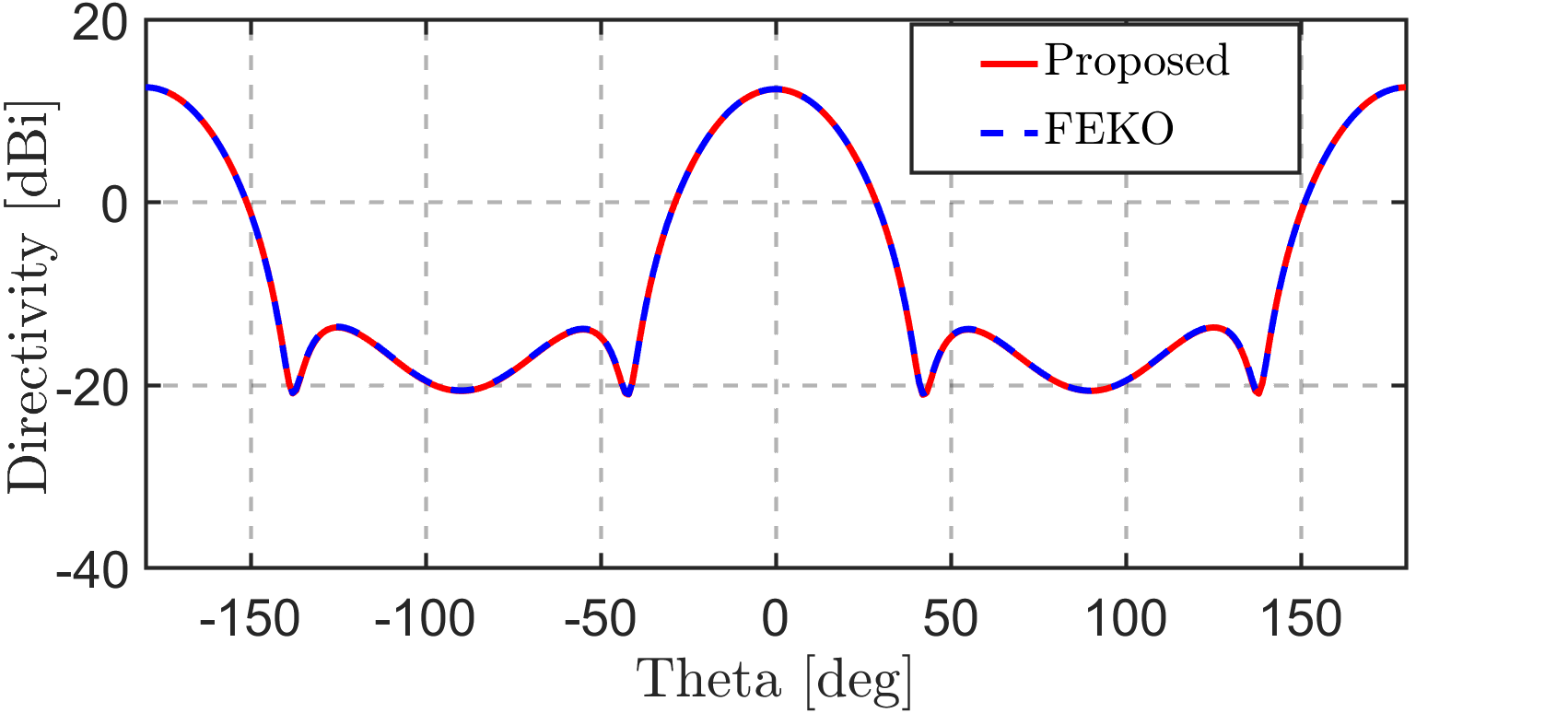}}} 
\\
\subfloat[\label{fig:Directivity_ex1_phi_0} ] 
{\resizebox{0.9\columnwidth} {!}
{\includegraphics{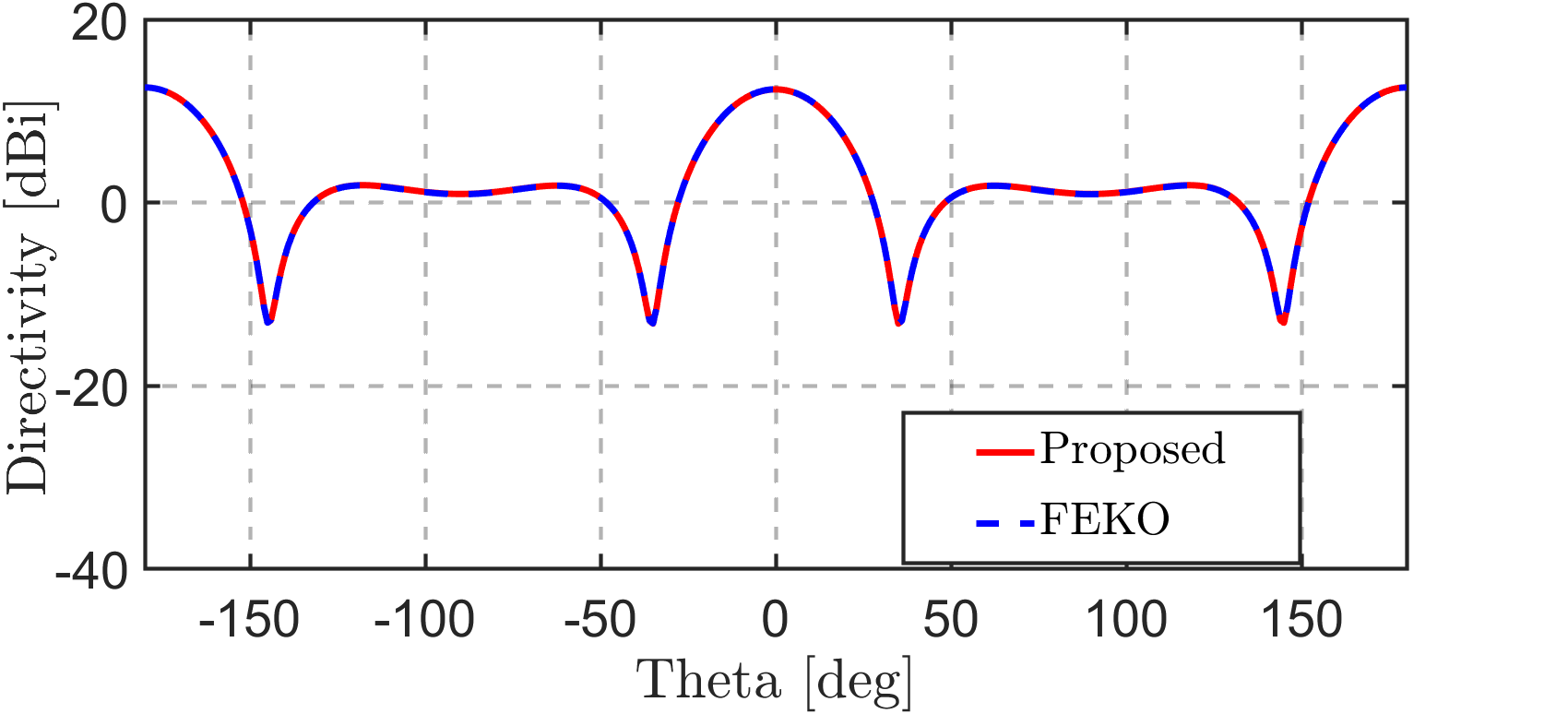}}} 
\end{center}
\caption{Directivity of the array of rectangular PEC patches calculated with the proposed method and FEKO~\cite{FEKO} for the (a) $\phi = 0^{\circ}$, (b) $\phi = 30^{\circ}$, and (c) $\phi = 90^{\circ}$ cuts.}
\label{fig:Directivity_example1}
\end{figure}

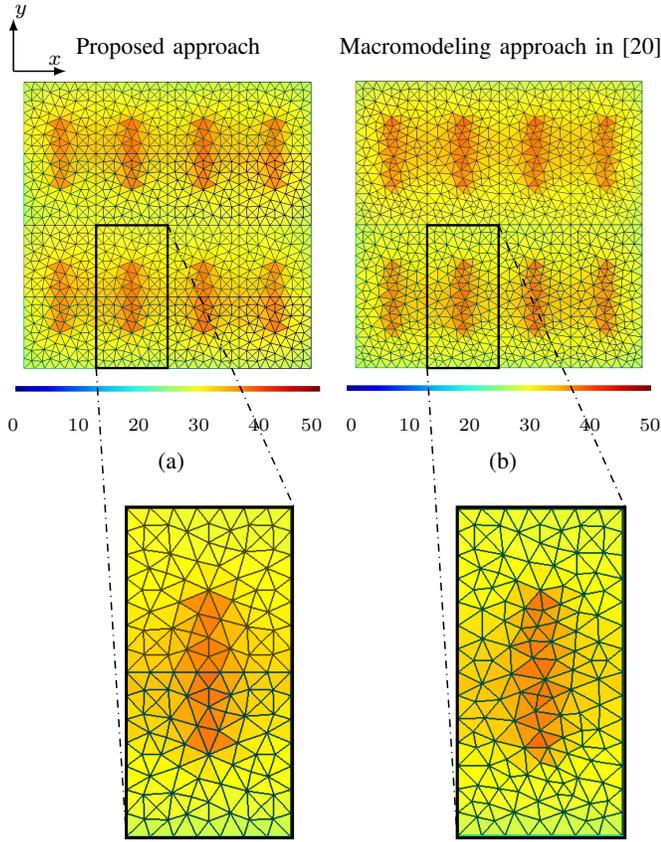
\begin{figure}[h]
\centering
\begin{tikzpicture}[scale=1]
\centering

\draw[anchor=west] node at (-1.2,-3.0) {\small{(a)}};
\draw[anchor=west] node at (3.2,-3.0) {\small{(b)}};
\draw[anchor=west] node at (-2.4,2.5) {\small{ Proposed approach}};
\draw[anchor=west] node at (1.1,2.5) {\small{ Macromodeling approach in \cite{Utkarsh_macromodel_2020}}};

\begin{scope}[shift={(-0.2,0.0)}]
\draw node at (-2.7,2.2+.8) {\footnotesize{$y$}};
\draw [line width=0.2mm,-latex] (-2.8,2.2) -- (-2.8,2.2+.7);
\draw node at (-2.8+.55,2.35) {\footnotesize{$x$}};
\draw [line width=0.2mm,-latex] (-2.8,2.2) -- (-2.8+.7,2.2);
\end{scope}

\node[anchor=west,inner sep=0] (image) at (-3.0,0.0){\includegraphics[width=0.225\textwidth]{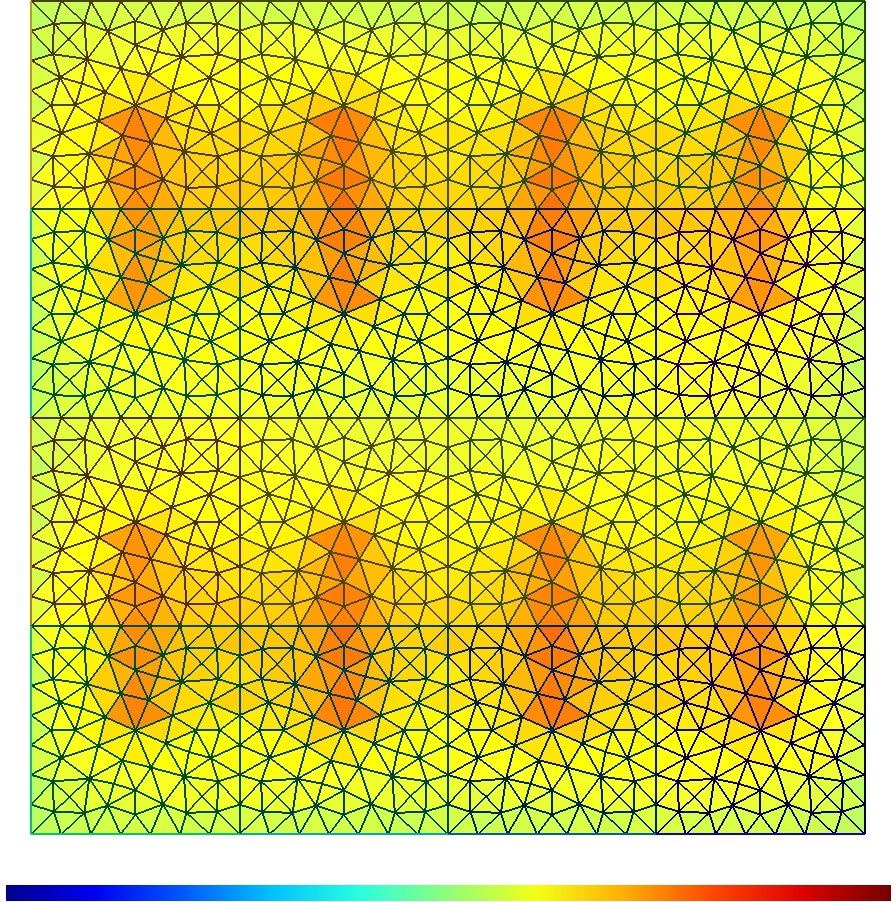}};
\draw[anchor=west] node at (-3.2,-2.5) {\scriptsize{$0$}};
\draw[anchor=west] node at (-3.2+4.0*1/5,-2.5) {\scriptsize{$10$}};
\draw[anchor=west] node at (-3.2+4.0*2/5,-2.5) {\scriptsize{$20$}};
\draw[anchor=west] node at (-3.2+4.0*3/5,-2.5) {\scriptsize{$30$}};
\draw[anchor=west] node at (-3.2+4.0*4/5,-2.5) {\scriptsize{$40$}};
\draw[anchor=west] node at (-3.2+3.9,-2.5) {\scriptsize{$50$}};

\draw[line width=0.3mm, rotate=0] {(-1.9,-1.75) rectangle ++(0.95,1.9)};
\node[anchor=west,inner sep=0] (image) at (-1.5,-5.8)
{\includegraphics[width=0.12\textwidth]{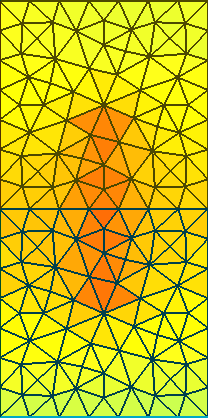}};
\draw[line width=0.4mm] (-1.5,-8.0) rectangle ++(2.21,4.4);
\draw [ dash dot, line width=0.2mm] (-1.9, -1.75) -- (-1.5,-8.0);
\draw [ dash dot, line width=0.2mm] (-1.9+0.95, -1.75+1.9) -- (-1.5+2.23,-8.0+4.4);

\begin{scope}[shift={(4.4,0.0)}]
\node[anchor=west,inner sep=0] (image) at (-3.0,0.0){\includegraphics[width=0.225\textwidth]{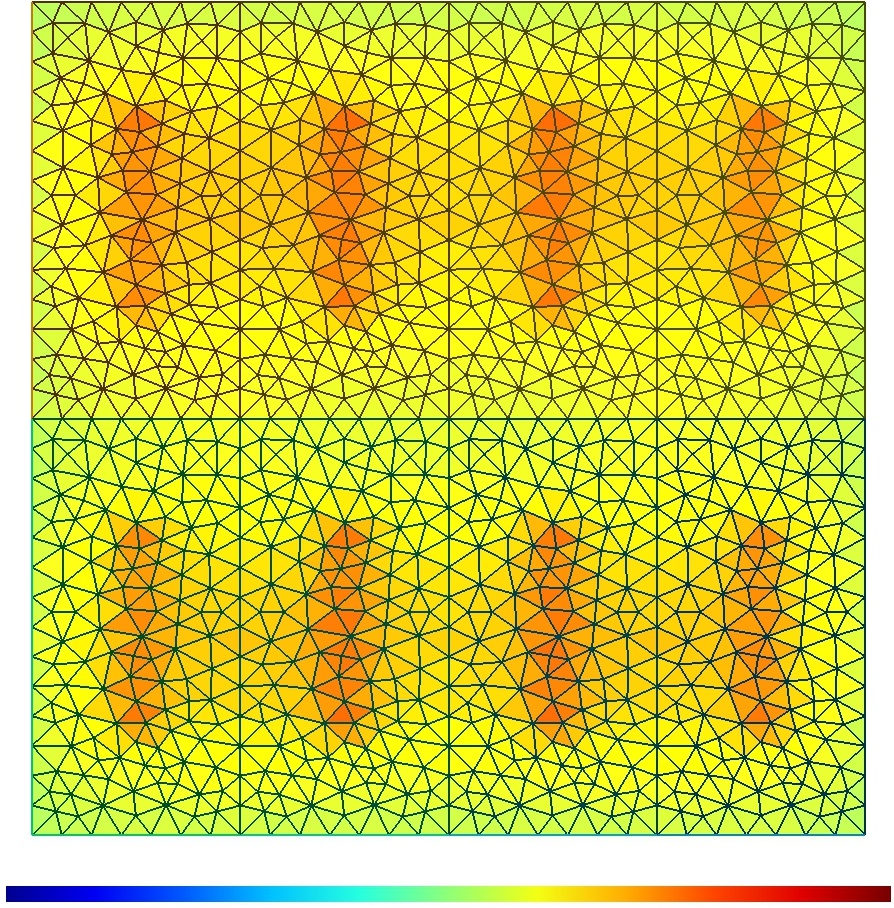}};
\draw[anchor=west] node at (-3.1,-2.5) {\scriptsize{$0$}};
\draw[anchor=west] node at (-3.2+4.0*1/5,-2.5) {\scriptsize{$10$}};
\draw[anchor=west] node at (-3.2+4.0*2/5,-2.5) {\scriptsize{$20$}};
\draw[anchor=west] node at (-3.2+4.0*3/5,-2.5) {\scriptsize{$30$}};
\draw[anchor=west] node at (-3.2+4.0*4/5,-2.5) {\scriptsize{$40$}};
\draw[anchor=west] node at (-3.2+4.0,-2.5) {\scriptsize{$50$}};

\draw[line width=0.3mm, rotate=0] {(-1.9,-1.75) rectangle ++(0.95,1.9)};
\node[anchor=west,inner sep=0] (image) at (-1.5,-5.8)
{\includegraphics[width=0.121\textwidth]{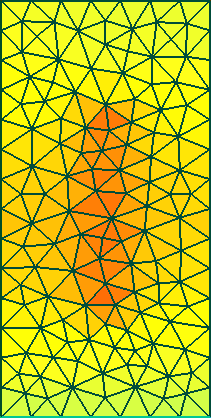}};
\draw[line width=0.4mm] (-1.5,-8.0) rectangle ++(2.23,4.4);
\draw [ dash dot, line width=0.2mm] (-1.9, -1.75) -- (-1.5,-8.0);
\draw [ dash dot, line width=0.2mm] (-1.9+0.95, -1.75+1.9) -- (-1.5+2.23,-8.0+4.4);
\end{scope}

\end{tikzpicture}
\caption{(a) Top view of the distribution of the equivalent electric current density (in $m\text{A}/\text{m}^2$) on the equivalent surfaces of (a) an array of $4 \times 4$ unit cells solved by the proposed approach (problem set of Fig.~\ref{fig:Array_Rectangular_Patches}) (b) an array of $4 \times 2$ unit cells solved by the previous macromodeling approach~\cite{Utkarsh_macromodel_2020} (problem set of Fig.~\ref{fig:Array_Rectangular_Patches2}). Part of the figure is magnified for better comparison of the current distribution between the two techniques.}
\label{fig:currentDistribution_Example1}
\end{figure}

The directivity of the array of rectangular PEC patches for the $\phi= 0^{\circ}$, $\phi= 30^{\circ}$, and $\phi= 90^{\circ}$ cuts is computed with the proposed solver. The results are plotted in Fig.~\ref{fig:Directivity_example1} and compared with FEKO-MLFMM. The two results match very well, validating the proposed approach. Also, distribution of the magnitude of the equivalent electric current density on the fictitious surfaces is depicted in Fig.~\ref{fig:currentDistribution_Example1} using the proposed solver and the result is compared to the one obtained with the previous macromodeling approach~\cite{Utkarsh_macromodel_2020}. The excellent agreement between the two approaches validates the accuracy of the current modeling in the proposed solver. 

\subsection{Three-Layer Linear-to-Circular Polarizer}
\label{Three-Layer Linear-to-Circular Polarizer}
\begin{figure}
\begin{center}
\scalebox{0.75}{
\includegraphics[width=.8\linewidth, angle=0]{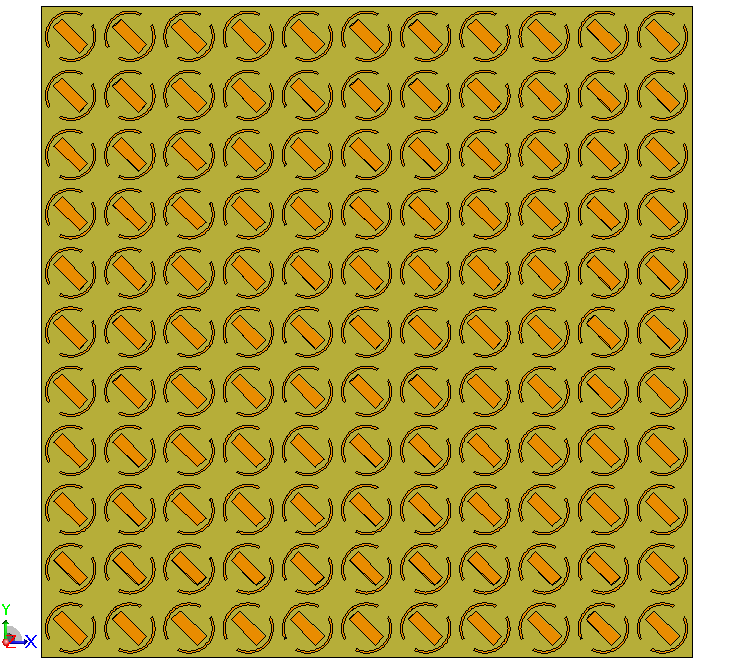} }
\end{center}
\caption{Top view of $11 \times 11$ polarizer.}
\label{fig:Polarizer_TopView}
\end{figure}

\begin{figure}
\begin{center}
\scalebox{0.85}{
\includegraphics[width=1.\linewidth, angle=0]{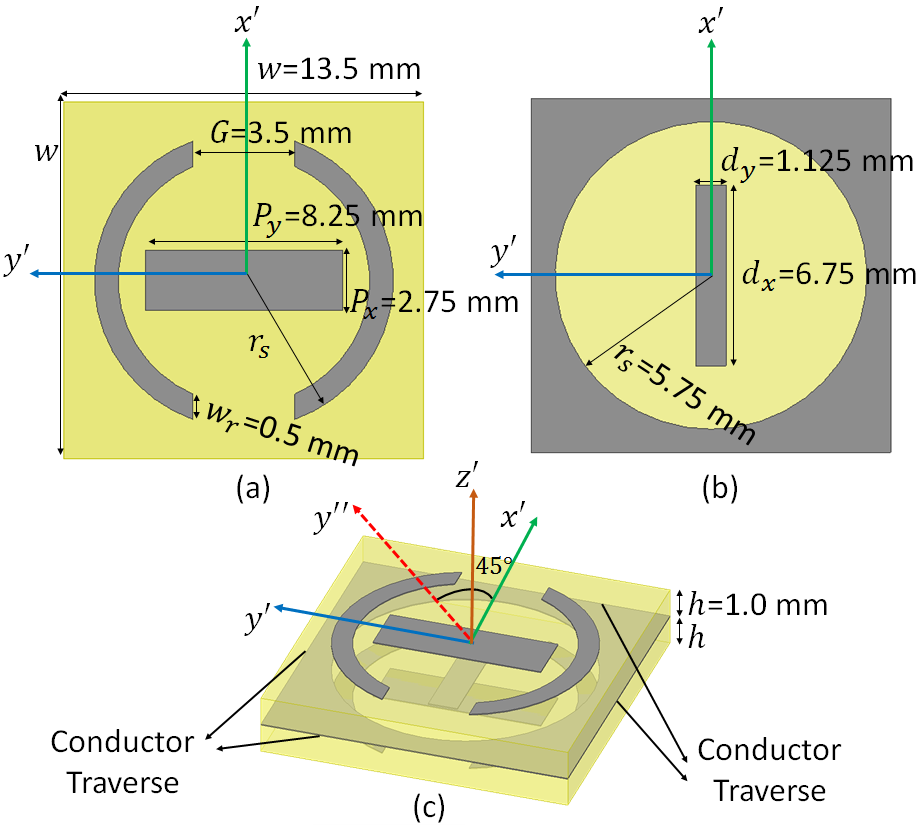} }
\begin{tikzpicture}[scale=1]
\end{tikzpicture}
\end{center}
\caption{Unit cell of the polarizer composed of three PEC layers stacked with two substrate layers. The PEC scatterers are shown with gray color. (a) The PEC scatterers of the top layer and the bottom layer are identical. (b) The middle layer, and (c) 3-D view of the unit cell which shows the PEC scatterers in the middle layer traverse the equivalent surface from four sides of the unit cell. }
\label{fig:Polarizer_unitcell}
\end{figure}

For the second example, we consider an $11 \times 11$ three-layer polarizer involving an inductive middle layer~\cite{EMS_Connected5}. The top view of this polarizer is shown in Fig.~\ref{fig:Polarizer_TopView}. The unit cell of the polarizer is shown in Fig.~\ref{fig:Polarizer_unitcell}. It is composed of three PEC layers stacked with two substrate layers, each with thickness of $1.0$ mm. Each unit cell has dimensions of $13.5$ mm $\times~13.5$ mm $\times~2.0$ mm. In order to convert linear polarization to circular polarization, the unit cell is rotated $45^{\circ}$ around its local $z'$-axis with respect to the main coordinate $(x',y',z')$ coordinate system, shown in Fig.~\ref{fig:Polarizer_unitcell}. Each unit cell is enclosed by an equivalent surface of size $13.5$ mm $\times~13.5$ mm $\times~4.0$ mm. The relative permitivity of each substrate layer is $\varepsilon_{r} = 2.2$. In this example, the polarizer is $4.75\lambda_0 \times 4.75\lambda_0$ at $9.6$ GHz, where $\lambda_0$ is the wavelength in free space. The structure is excited by the same dipole used for the array of rectangular PEC patches in Section~\ref{Example1_Array_Rectangular_Patches}.  

As shown in Fig.~\ref{fig:Polarizer_unitcell}, the PEC traces of the complementary structure traverse the equivalent surface from four sides which makes the simulation of this structure difficult for EPA approach in SIE, due to the need for accurate current modeling of the PEC traces cut by the equivalent surfaces. Furthermore, the solution of such multiscale arrays becomes prohibitively expensive for traditional SIE methods, such as the Poggio-Miller-Chang-Harrington-Wu-Tsai (PMCHWT) formulation~\cite{Gibson-MoM}, in terms of required CPU time and memory.  

\begin{figure}[t]
\begin{center}
\subfloat[\label{fig:Directivity_ex2_phi_0} ] 
{\resizebox{0.8\columnwidth} {!}
{\includegraphics{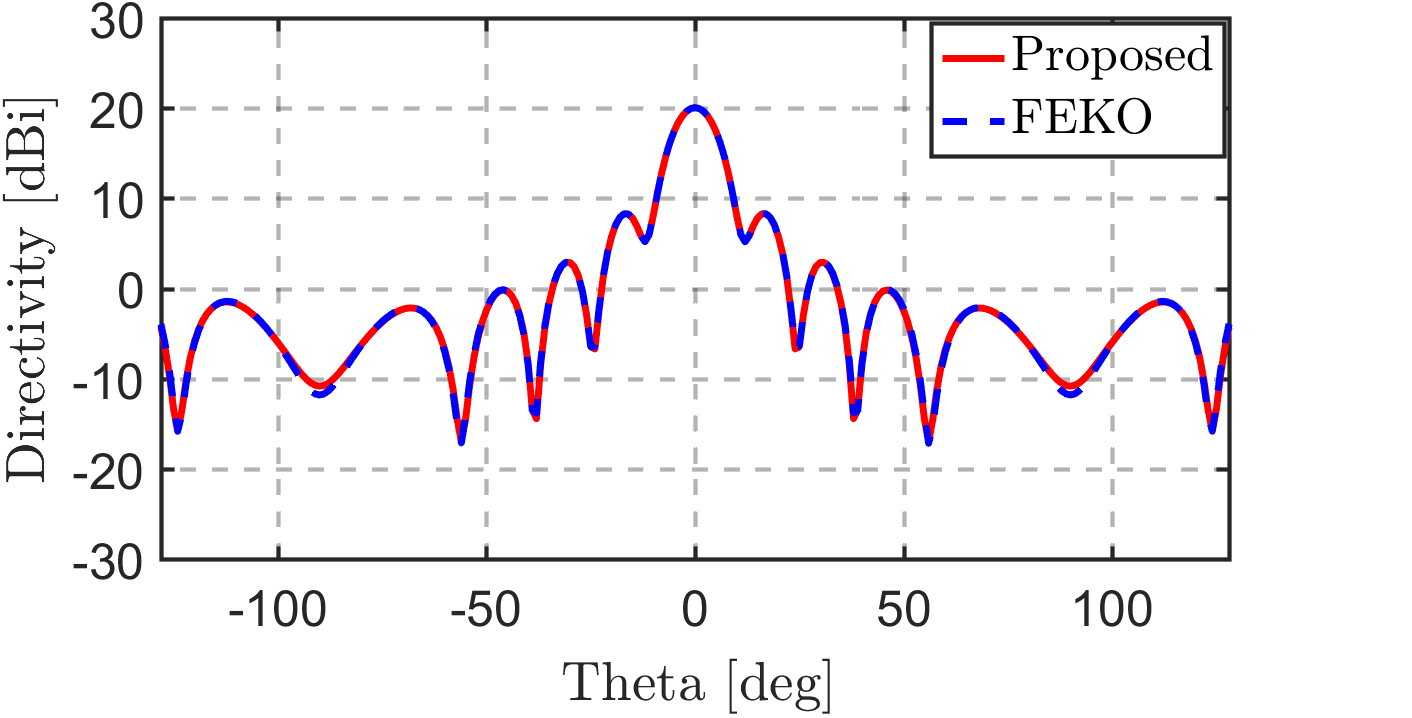}}} 
\\
\subfloat[\label{fig:Directivity_ex2_RHCP} ] 
{\resizebox{0.8\columnwidth} {!}
{\includegraphics{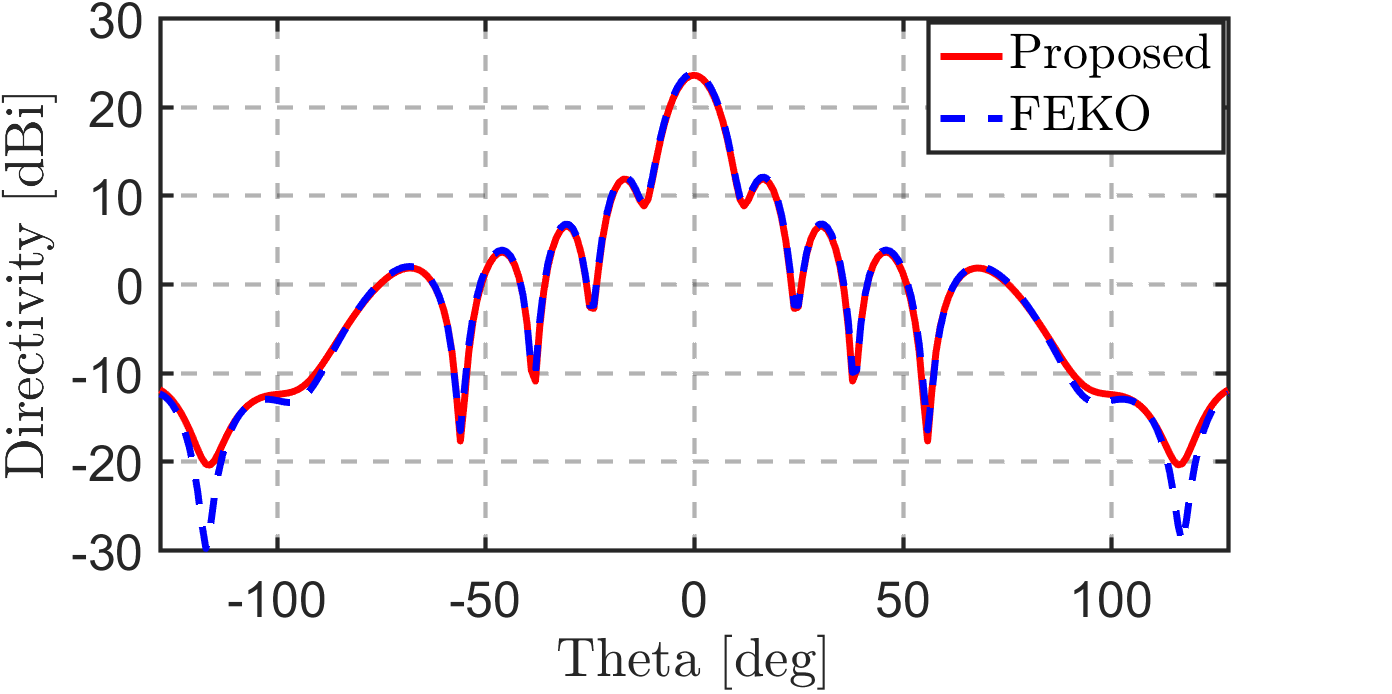}}} 
\\
\subfloat[\label{fig:Directivity_ex2_LHCP} ] 
{\resizebox{0.8\columnwidth} {!}
{\includegraphics{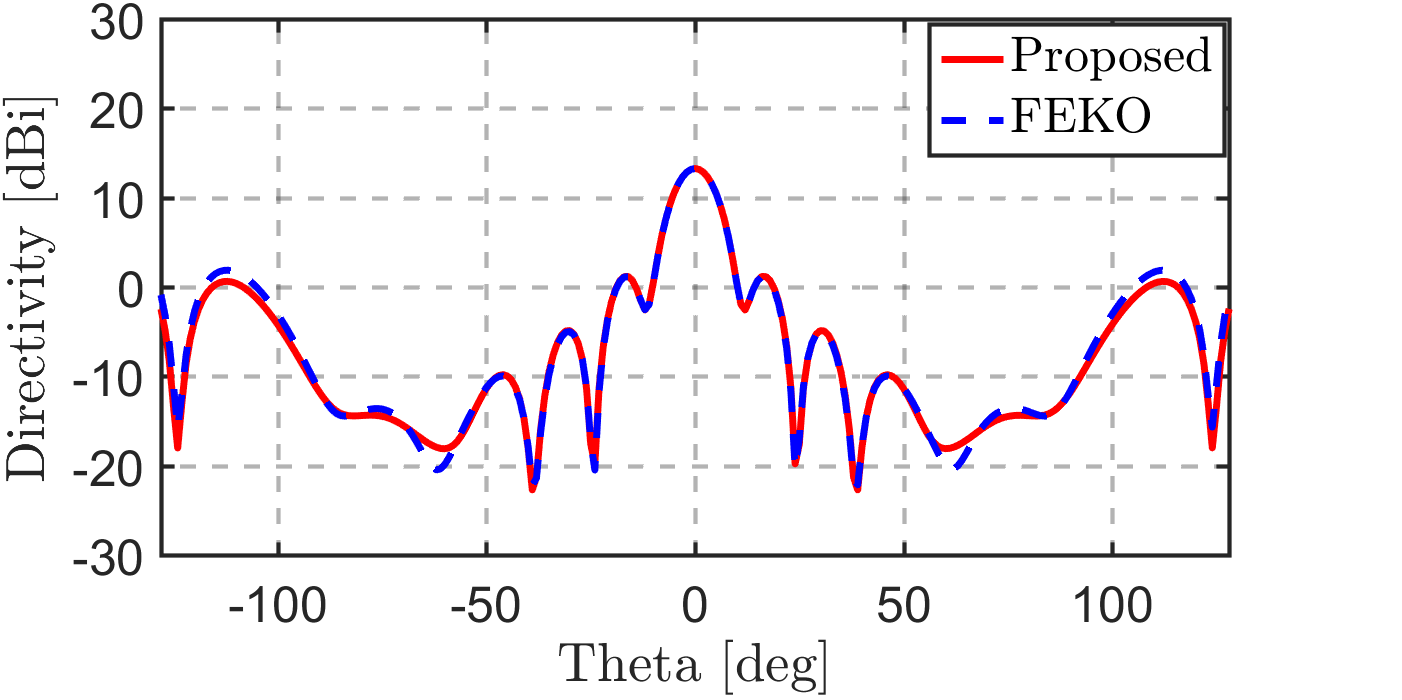}}} 
\end{center}
\caption{Directivity of the $11 \times 11$ polarizer calculated with the proposed method and FEKO~\cite{FEKO} for the $\phi = 0^{\circ}$ cut. (a) Total field, (b) RHCP pattern, and (c) LHCP pattern.}
\label{fig:Directivity_polarizer}
\end{figure}

In our simulation, the PEC scatterers of the unit cell are discretized with triangular meshes with a characteristic length of $0.95$ mm while a characteristic length of $2$ mm is chosen for the dielectric regions. Also, the equivalent surface in the macromodeling approach is discretized with a characteristic length of $2$ mm. We calculated the scattered field from the polarizer using the proposed solver and FEKO-MLFMM. The final design contains $174,380$ and $449,928$ unknowns for the proposed approach and FEKO simulation, respectively. The directivity of the three-layer polarizer in the $\phi= 0^{\circ}$ cut is plotted in Fig.~\ref{fig:Directivity_polarizer}. Fig.~\ref{fig:Directivity_polarizer}\subref{fig:Directivity_ex2_phi_0} shows the radiation pattern of the total electric field. Fig.~\ref{fig:Directivity_polarizer}\subref{fig:Directivity_ex2_RHCP} and Fig.~\ref{fig:Directivity_polarizer}\subref{fig:Directivity_ex2_LHCP} show the right-hand circular polarization (RHCP) pattern obtained using $\gv{E}_{\text{RHCP}}={(\gv{E}_{\phi}- j\gv{E}_{\theta})}{/\sqrt{2}}$ and the left-hand circular polarization (LHCP) patterns obtained using $\gv{E}_{\text{LHCP}}={(\gv{E}_{\phi}+ j\gv{E}_{\theta})}{/\sqrt{2}}$, respectively. An excellent match between the proposed solver and the FEKO commercial solver validates the accuracy of the proposed method. The memory consumption and timing results of both methods are compared in Table~\ref{tab:polarizer}. As summarized in the Table~\ref{tab:polarizer}, the proposed macromodeling approach is $28$ times faster and requires $5$ times less memory than FEKO.    

\begin{table}[t]
\caption{Simulation statistics for the $11 \times 11$ polarizer considered in Section~\ref{Three-Layer Linear-to-Circular Polarizer}}
\label{tab:polarizer}
\begin{center}
\begin{tabular}{|l| c | c|}
\hline
& FEKO-MLFMM & Proposed \\
\hline
\multicolumn{3}{|c|}{Memory Consumption}\\
\hline
Total number of unknowns & 449,928 & 174,380 \\
Memory used & 316.9~GB & 64.27~GB \\
\hline
\multicolumn{3}{|c|}{Timing Results}\\
\hline
Macromodel generation & N/A & 3~min \\
Matrix fill time & 1.56~h &  19~min \\
Preconditioner factorization & 23.1~h & 34~min \\
Iterative solver & 6.3~h & 9~min \\
Total computation time & 31.9~h & 1.10~h \\
\hline
\end{tabular}
\end{center}
\end{table}

\subsection{Phase-Rotation Transmitarray}
\label{Phase-Rotation Transmitarray}
\begin{figure}
\begin{center}
\scalebox{0.65}{
\includegraphics[width=1.2\linewidth, angle=0]{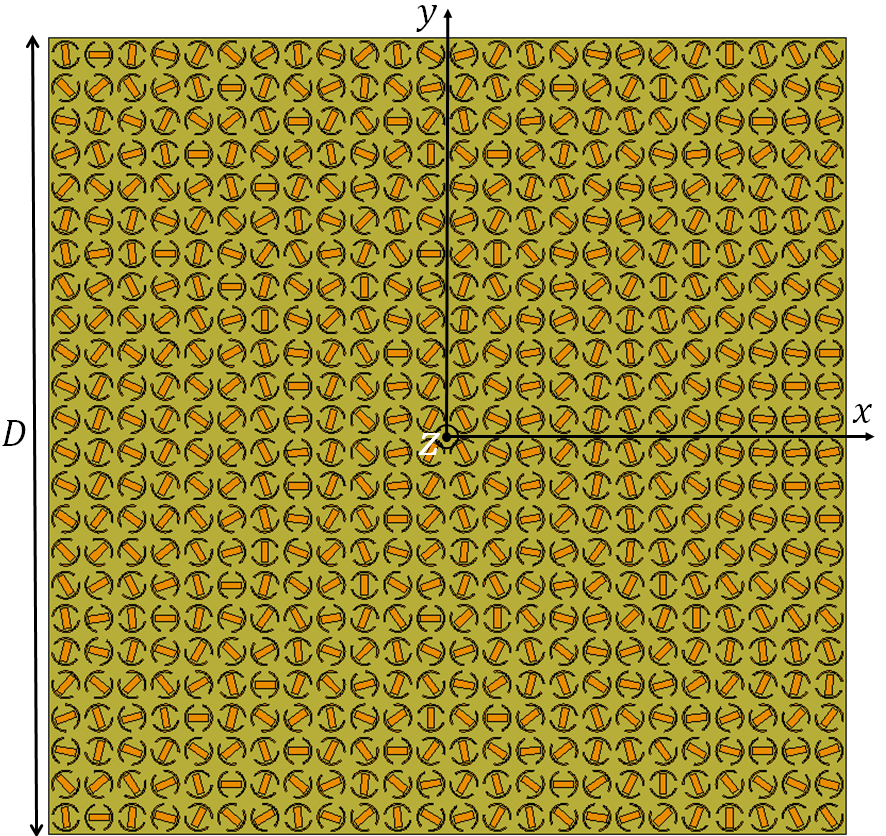} }
\end{center}
\caption{Top view of $24 \times 24$ phase rotation TA with $D=120$ mm. The structure is symmetric in relation to the $x$-axis and is designed to radiate the main beam in $\alpha_0 = 32.5^{\circ}$ (with respect to the $z$-axis) direction. The TA has $36$ unique unit cells.}
\label{fig:TA_TopView}
\end{figure}

\begin{figure}
\begin{center}
\scalebox{0.85}{
\includegraphics[width=1.\linewidth, angle=0]{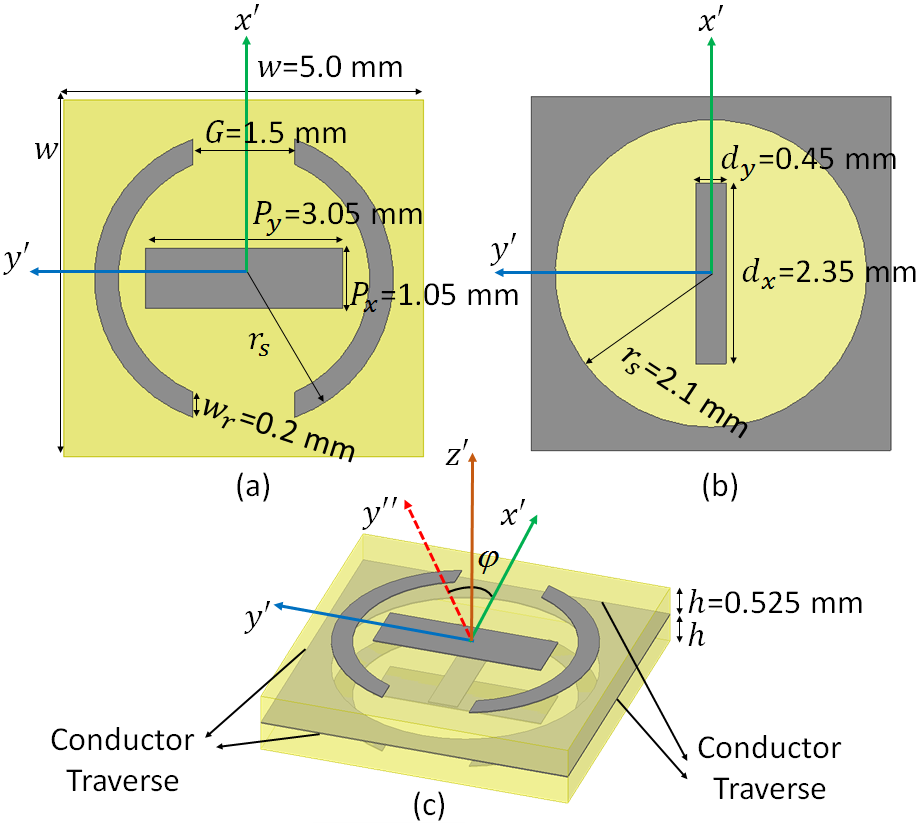} }
\begin{tikzpicture}[scale=1]
\end{tikzpicture}
\end{center}
\caption{Unit cell of the phase rotation TA. (a) Geometry of the PEC scatterers in the top and bottom layer, which are identical, (b) geometry of the PEC scatterers in the middle layer, and (c) 3-D view of the entire unit cell. The rotation angle of each unit cell is $\varphi$ and obtained by \eqref{rotation_angle}.}
\label{fig:TA_unitcell}
\end{figure}

Here, we consider a $24 \times 24$ transmitarray (TA) with unit cells implementing phase shifts using sequential rotation~\cite{EMS_Connected8}. This example was chosen to demonstrate the performance of the proposed macromodeling approach for the simulation of EM surfaces with a heterogeneous collection of unit cells composing its surface. The top view of the phase rotation TA, centered at the origin, is shown in Fig.~\ref{fig:TA_TopView}. The TA substrate has two layers, each with a thickness of $0.525$ mm. The structure is symmetric with respect to the $x$-axis. The unit cell of the TA is shown in Fig.~\ref{fig:TA_unitcell}. Each unit cell has dimensions of $5.0$ mm $\times~5.0$ mm $\times~1.05$ mm. The relative permittivity of the bottom and top layers of the substrate is $\varepsilon_r = 2.2$. In this example, the TA is electrically large: $12\lambda_0 \times 12\lambda_0$ at $30$ GHz, where $\lambda_0$ is the wavelength in free space. The model is excited by a horn antenna. 
The horn antenna is placed $0.1$ m away from the TA along $z$-axis with $F/D = 0.83$ (see Fig.\ref{fig:TA_TopView}). 

As shown in Fig.~\ref{fig:TA_unitcell}, each unit cell is rotated about its local $z'$-axis by an angle
\begin{equation}
\label{rotation_angle}
\begin{aligned}
\varphi_{\text{unit cell}}(x,y) = -\frac{1}{2}k_0 \left[ x \sin(\alpha_0) - \sqrt{x^2+y^2+F^2} \right],
\end{aligned}
\end{equation}
where $(x,y)$ is the centroid of the unit cell and $\alpha_0 = 32.5^{\circ}$ is the desired beam position. In~\eqref{rotation_angle}, $F=0.1$ m is the distance between the phase center of the spherical waves from the horn antenna and the surface of the TA. The rotation angles calculated with~\eqref{rotation_angle} were rounded to multiples of 5 degrees. This level of discretization is sufficient to achieve a proper design~\cite{EMS_Connected8}. This way, the array in Fig.~\ref{fig:TA_TopView} is composed of $36$ unique unit cells. 

To calculate the RHCP and LHCP patterns of the phase rotation TA, it needs to be excited by both $x$- and $y$-polarized incident fields. The CP patterns can then be obtained using
\begin{equation}
\label{E_RHCP}
\begin{aligned}
E_{\text{RHCP}} = \frac{E_x^{\theta}-E_y^{\phi}+ j\left(E_x^{\phi}+E_y^{\theta}\right)}{2} 
\end{aligned}
\end{equation}
and 
\begin{equation}
\label{E_LHCP}
\begin{aligned}
E_{\text{LHCP}} = \frac{E_x^{\theta}+E_y^{\phi}+ j\left(-E_x^{\phi}+E_y^{\theta}\right)}{2}, 
\end{aligned}
\end{equation}
respectively. In \eqref{E_RHCP} and \eqref{E_LHCP}, $E_a^b$ is the pattern of the $E$-field in the $b$-direction ($\phi$ or $\theta$) generated by an $x$- or $y$-polarized horn antenna. 

\begin{figure}[t]
\begin{center}
\subfloat[\label{fig:Directivity_ex3_RHCP} ] 
{\resizebox{0.85\columnwidth} {!}
{\includegraphics{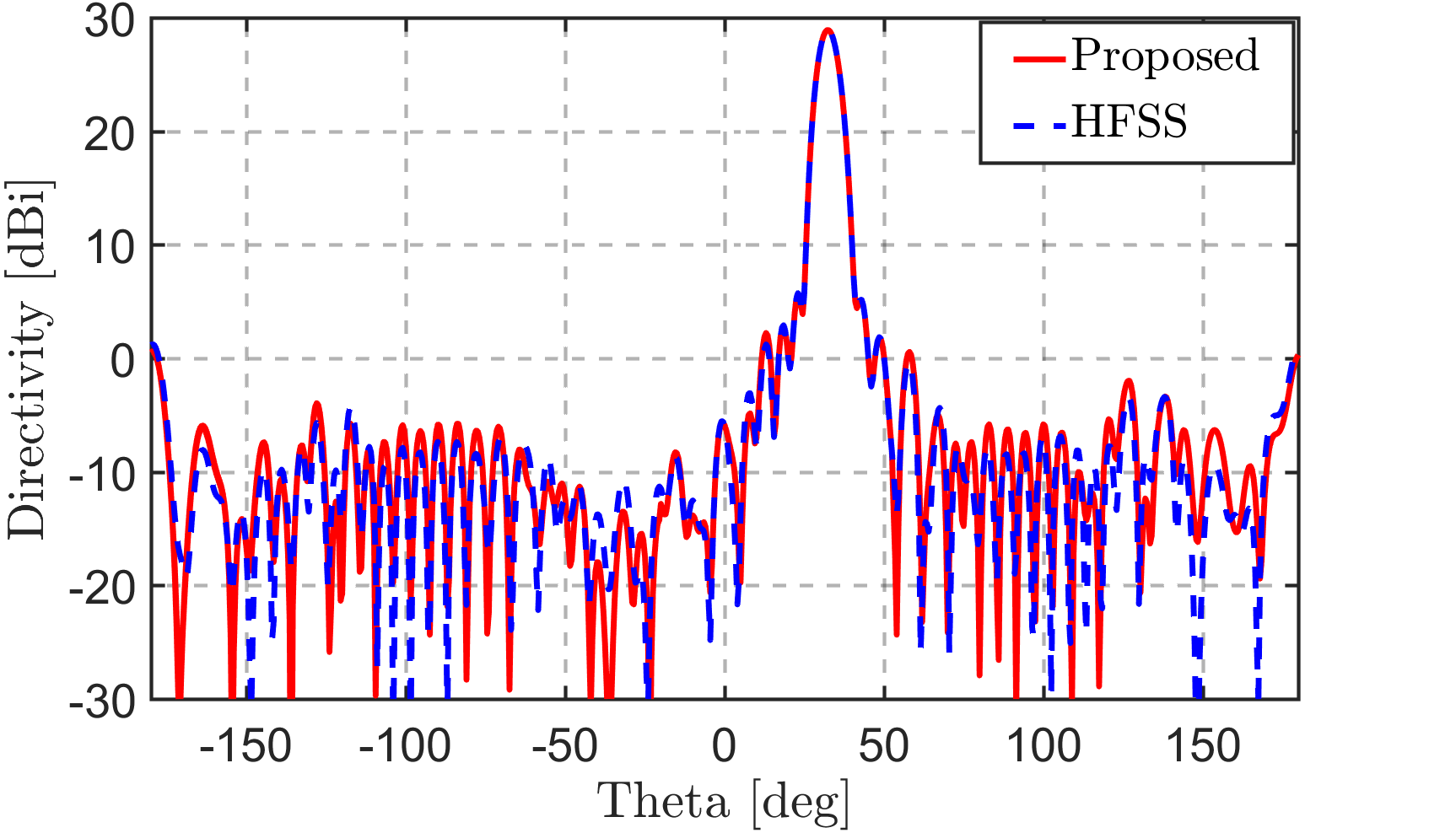}}} 
\\
\subfloat[\label{fig:Directivity_ex3_LHCP} ] 
{\resizebox{0.85\columnwidth} {!}
{\includegraphics{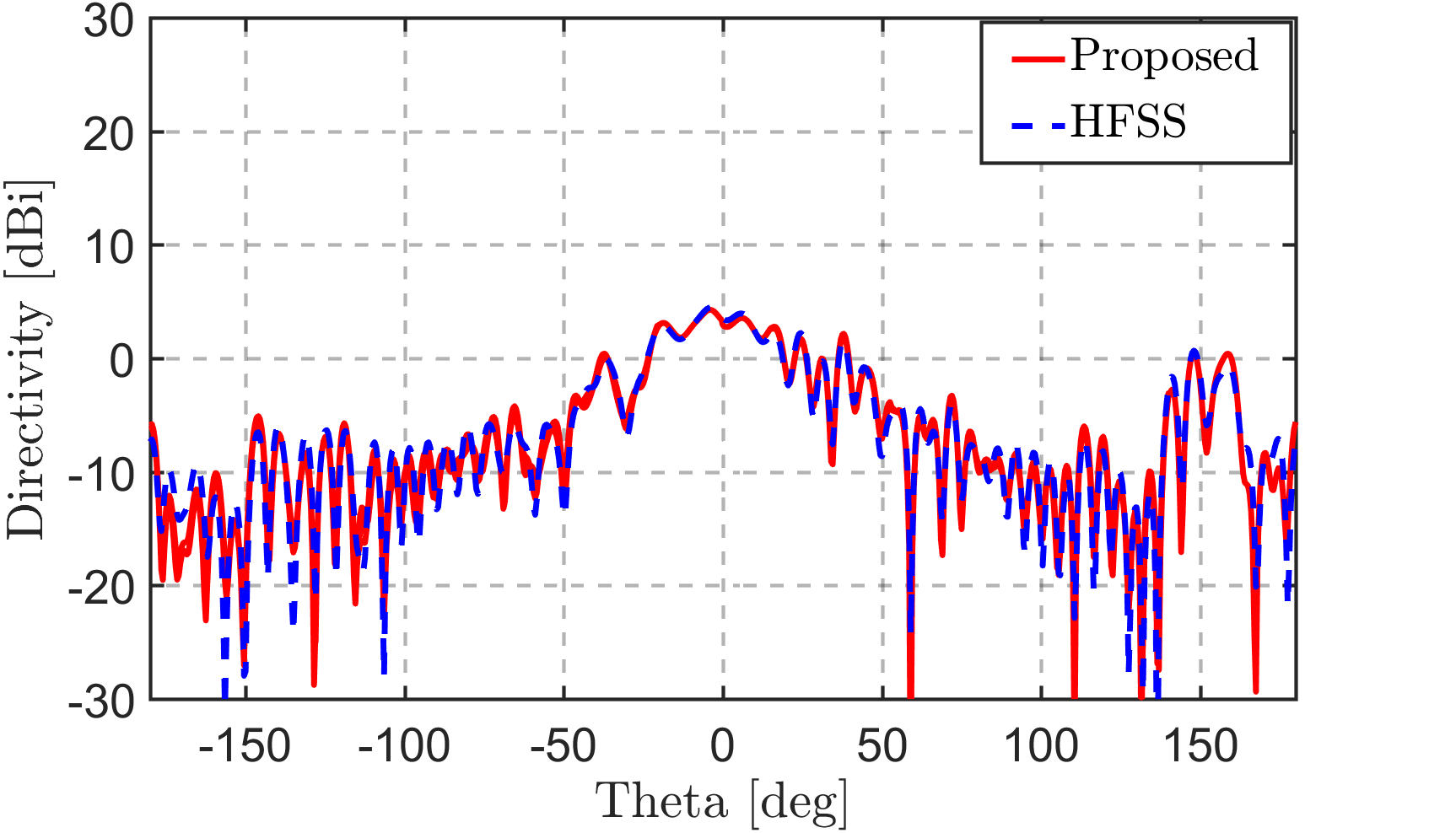}}} 
\end{center}
\caption{Directivity of the $24 \times 24$ phase rotation TA calculated with the proposed method and HFSS in the $\phi = 0^{\circ}$ cut. (a) RHCP pattern and (b) LHCP pattern.}
\label{fig:TA_result}
\end{figure}

In the proposed macromodeling approach, the PEC scatterers inside each unit cell are discretized with a characteristic length of $0.5$ mm, while a characteristic length of $1$ mm is chosen for the dielectric regions and the equivalent surfaces. Hence, this discretization produces $1,368$ unknowns for each unit cell, which leads to $506,888$ total unknowns. To simulate this structure using a commercial solver, the hybrid finite element-boundary integral (FE-BI) technique is employed in Ansys HFSS~\cite{HFSS}. The RHCP and LHCP radiation patterns of the TA for the $\phi = 0^{\circ}$ cut obtained from the proposed solver are plotted in Fig.~\ref{fig:TA_result}\subref{fig:Directivity_ex3_RHCP} and Fig.~\ref{fig:TA_result}\subref{fig:Directivity_ex3_LHCP}, respectively, and the results are compared with HFSS results. An excellent match between the proposed solver and HFSS FE-BI solver validates the accuracy of the proposed method. The memory and time requirements of both solvers are listed in Table~\ref{tab:phaserotation TA}. For this simulation, the proposed macromodeling technique required $172.3$ GB memory, while the HFSS solver required $174.8$ GB memory. Also, the proposed approach took $3.45$ h to simulate this structure, while the HFSS took $25.42$ h. Hence, the proposed macromodeling approach is $7$ times faster than HFSS solver while requiring the same amount of memory.  

\begin{table}[t]
\caption{Simulation statistics for the $24 \times 24$ phase rotation TA considered in Section~\ref{Phase-Rotation Transmitarray}}
\label{tab:phaserotation TA}
\begin{center}
\begin{tabular}{|l| c | c|}
\hline
& HFSS & Proposed \\
\hline
\multicolumn{3}{|c|}{Memory Consumption}\\
\hline
Total number of unknowns & $-$ & 506,888 \\
Memory used & 174.8~GB & 172.3~GB \\
\hline
\multicolumn{3}{|c|}{Timing Results}\\
\hline
Macromodel generation & $-$ & 44~min \\
Matrix fill time & $-$ &  45~min \\
Preconditioner factorization & $-$ & 46~min \\
Iterative solver & $-$ & 55~min \\
Total computation time & 25.42~h & 3.45~h \\
\hline
\end{tabular}
\end{center}
\end{table} 

\section{Conclusion}
\label{Conclusion}
The macromodeling approach is generalized for the simulation of EM surfaces composed of unit cells connected via PEC traces. We presented a macromodel, through which a complex unit cell can be replaced by equivalent electric and magnetic surface current densities on a simple fictitious surface enclosing the unit cell. In the proposed macromodeling approach, the PEC traces of a unit cell are allowed to traverse the fictitious surface, which enables the efficient analysis of different EM surfaces. In particular, we demonstrated how to tackle problems where a current-carrying conductor of a unit cell traverses the equivalent surfaces. The new method based on half RWG basis functions accurately models the continuity of the surface currents and can be extended to structures with different conductor traverse positions. 

The proposed approach generates identical basis functions on the equivalent surfaces of particular EM surfaces. This helps us to represent the interaction matrix of uniform array of equivalent surfaces with a Toeplitz form even when unit cells are not the same. This property allows us to accelerate the simulation of large EM surfaces via FFT, while still rigorously capturing the mutual coupling between the unit cells. The proposed method is shown to produce accurate results while reducing the time and resources needed to analyze large EM surfaces with respect to commercial solvers. The proposed acceleration scheme based on the FFT requires some regularity in the way conductors traverse the boundaries between unit cells. Future work will remove this limitation, devising an alternative acceleration scheme valid for conductors traversing equivalent surfaces in arbitrary positions. 



\bibliographystyle{_template/IEEEtran}
\bibliography{_template/IEEEabrv,_template/IEEEfull,library_fixed}

\end{document}